\documentclass{siamltex}
\title{Coarse-graining schemes and  a posteriori error estimates
for stochastic lattice systems}

\author{Markos A. Katsoulakis\thanks{Department of Mathematics,
University of Massachusetts ({\tt markos@math.umass.edu}).}
        \and Petr Plech\'a\v{c}\thanks{Mathematics Institute,
University of Warwick ({\tt plechac@maths.warwick.ac.uk}).}
        \and Luc Rey-Bellet\thanks{Department of Mathematics,
University of Massachusetts({\tt lr7q@math.umass.edu}).}
        \and Dimitrios K. Tsagkarogiannis\thanks{
Max Planck Institute for Mathematics in the Sciences,
({\tt tsagkaro@mis.mpg.de}).}
}

\usepackage{amsfonts,amssymb}        
\usepackage{amsmath}
\usepackage{epsfig}

\def\RANGE{L}                            
\def\RANG_S{S}                            
\def\CGPARQ{Q}                            
\def\BARIT#1{{\bar {#1}}}
\def\BARH{\BARIT{H}}
\def\LATT{{\Lambda}_N}
\def\SIGMA{{\mathcal{S}_N}}
\def\LATTC{{\bar{\Lambda}_{M}}}
\def\SIGMAC{\bar{\mathcal{S}}_{M}} 

\def\alfa{\alpha}
\def\vita{\omega}

\def\SPINSPC{\Sigma^c}

\def\PROCMICRO{\sigma_t}

\def\PROCMIC{\{\sigma_t\}_{t\geq 0}}
\def\PROCMAC{\{\eta_t\}_{t\geq 0}}

\def\JNOT{V}

\def\STATESP{\mathcal{S}}
\def\CUBE{C}
\def\COP{\mathbf{F}}

\def\BARH{\bar H}

\def\RELENTR{\mathcal{R}}
\def\RELENT#1#2{\mathcal{R}\left({#1}\SEP{#2}\right)}

\def\Rr{{\mathcal{R}}}
\def\Gg{\mathcal{G}}

\def\Oo{\mathcal{O}}

\def\Tt{\mathcal{T}}

\def\R{\mathbb{R}}

\def\T{\mathbb{T}}

\def\Z{\mathbb{Z}}

\def\DT{\Delta t}

\def\VEC{\mathbf}
\def\EXP#1{e^{#1}}

\def\EXPECT{{\mathbb{E}}}

\def\SPACE{\;\;\;}          
\def\COMMA{\,,}             
\def\PERIOD{\,.}            
\def\SEP{{\,|\,}}           
\def\SUBSECT#1{{\smallskip\noindent{\it{#1}\/}.}}
\def\VIZ#1{(\ref{#1})}      

\def\BIGO{\Oo}

\def\qed{$\Box$}
\def\PROOF{\noindent{\sc Proof:}\ }

\def\VEC#1{{\mathbf{#1}}}        

\def\NORM#1{\Vert\,#1\,\Vert}

\def\HALF{\frac{1}{2}}

\def\GRADV{\nabla V}

\def\LP{\left(}
\def\RP{\right)}

\newtheorem{remark}{Remark}[section]

\newtheorem{scheme}{Scheme}[section]

\def\supp{\mbox{supp\,}}
\def\card{\mbox{card\,}}
\def\BARJ{\bar J}

\begin{document}

\maketitle

\begin{abstract}
The primary objective  of  this work is to  develop coarse-graining
schemes for stochastic many-body microscopic models and quantify their
effectiveness in terms of a priori and a posteriori  error analysis.  In
this paper  we focus on  stochastic lattice systems of
interacting particles at equilibrium. 
The proposed algorithms are derived from  an initial  coarse-grained
approximation  that is directly computable by Monte Carlo simulations, 
and the corresponding numerical error
is calculated using the specific relative entropy between the exact
and approximate coarse-grained equilibrium measures.  
Subsequently we carry out a cluster expansion around this first--and 
often inadequate--approximation and
obtain more accurate coarse-graining schemes.
The cluster expansions yield also sharp a posteriori error estimates for
the coarse-grained approximations that can be used for the construction of
adaptive coarse-graining methods.
We present a number of numerical examples that demonstrate that the
coarse-graining schemes developed here allow for
accurate predictions of critical behavior and hysteresis  in systems with
intermediate and long-range interactions. We also present
examples where they substantially improve
predictions of earlier coarse-graining schemes for short-range
interactions.
\end{abstract}



\pagestyle{myheadings}
\thispagestyle{plain}
\markboth{M.A. Katsoulakis, P. Plech\'a\v{c}, L. Rey-Bellet, D. K. Tsagkarogiannis}%
         {Coarse-graining schemes and a posteriori estimates}
%
%
%
%
\section{Introduction}\label{intro}
In the recent years there has been a growing interest in developing
and analyzing hierarchical coarse-graining methods for the purpose of
modeling and simulation across scales for systems arising in a broad
spectrum of scientific disciplines ranging from materials science to
biology and atmosphere/ocean science.  Typically in microscopic
simulations of complex systems the model is formulated in terms of
simple rules describing interactions between small-scale degrees of
freedom such as individual particles or spin variables. On the other
hand computational difficulties arise immediately from evaluating
their interactions for any realistic size spatio-temporal scales.
When a coarse-grained model becomes available, it has fewer
observables than the original microscopic system making it
computationally more efficient than the direct numerical simulations.
At the same time it is expected that it can describe accurately the
unresolved degrees of freedom.  The coarse-graining strategy we are
pursuing here is attempting to address these goals in the context of
equilibrium sampling of stochastic lattice systems by combining
methods from statistical mechanics and
perturbation analysis.

In this paper we consider stochastic lattice systems such as
Ising-type models as a paradigm of hierarchical coarse-graining that
can provide an explicit numerical method with prescribed error
tolerance. Such lattice systems for $N$ particles are defined in terms
of a microscopic lattice Hamiltonian $H_N(\sigma)$ with $\sigma$ being
the microscopic configuration and an a-priori Bernoulli measure
$P_N(d\sigma)$.  We perform a coarse-graining by subdividing the
lattice into coarse cells and defining variables $\eta$ on each coarse
cell to be the total magnetization in the cell. The corresponding
renormalization group map (known as the Kadanoff transform) \cite{Golden,Kad} is 
defined by the formula
\[
  \EXP{-\beta \BARH_M(\eta)}=\int \EXP{-\beta H_N(\sigma)} P_N(d\sigma|{\eta})\COMMA
\]
where $\BARH_M(\eta)$ is the Hamiltonian at the coarse level and 
$P_N(d\sigma|\VEC{\eta})$ is the conditional probability of having 
a microscopic configuration $\sigma$ given a configuration $\eta$ at the 
coarse level. 

Such a Hamiltonian defined on a coarser level than the microscopic, is
an exact equivalent of the microscopic Hamiltonian $H_N$, in the sense
that the finer degrees of freedom have been averaged.  However, it cannot
be easily calculated explicitly and hence used in numerical
simulations. Our perspective is to approximate it by viewing it as a
perturbation of a coarse-grained approximating Hamiltonian
$\BARH_M^{(0)}$ suggested in \cite{KMVPNAS,KMV} and defined
by
\[
   \BARH_M^{(0)}(\eta) \,=\, \int H_N(\sigma) P_N(d\sigma|{\eta})\,.
\]
A closely related coarse-grained Hamiltonian was suggested
independently in \cite{steph, steph2}, where it was constructed in an
equilibrium context using a wavelet expansion.  

Using this first approximation we have 
\[
   \BARH_M(\eta) \,=\, \BARH_M^{(0)}(\eta) -\frac{1}{\beta} \log \int  
   \EXP{-\beta (H_N(\sigma)-\BARH_M^{(0)}(\eta))} P_N(d\sigma|{\eta})\PERIOD
\]
The fact that the conditional probability $P_N(d\sigma|{\eta})$ factorizes 
at the level of the coarse cells allows us to use cluster expansion 
techniques  to write a series expansion for $\BARH_M(\eta)$ around 
$\BARH_M^{(0)}$. We obtain the
following series
\begin{equation}\label{intro_expansion}
   \BARH_M(\eta)=\BARH_M^{(0)}(\eta)+\BARH_M^{(1)}(\eta)
               +\cdots+\BARH_M^{(p)}(\eta) + \BIGO(\epsilon^{p+1})\COMMA\;
               p=2,\dots
\end{equation}
uniformly in $\eta$, where the correction terms $\BARH_M^{(1)}(\eta)$,
$\BARH_M^{(2)}(\eta)$ etc. can be calculated explicitly with the
relevant errors and $\epsilon$ is a small parameter depending on the
characteristics of the coarse-graining, the potential and the
inverse temperature. In this paper we first show that this strategy works well
provided that the interactions have a range which is long compared to
the size of the coarse cells and that they vary slowly over the size
of a coarse cell.  The parameter $\epsilon$ 
(given explicitly in \VIZ{small_parameter})
encapsulates these conditions. 

Notice that we do {\em not} perform the
usual high-temperature expansion using the Bernoulli product measure
$P_N(d\sigma)$ but rather expand the Hamiltonian around a well-chosen
first approximation using the product structure at the coarse level of
the conditional probability $P_N(d\sigma|\VEC{\eta})$. This allows us
to construct these approximating Hamiltonians well-beyond the
temperature range allowed in a standard high-temperature
expansion. The basic tool for cluster expansions is the reformulation
of the system in terms of what is called the polymer model.  This
technique originates from Mayer \cite{Mayer} and Peierls
\cite{Peierls} for the case of high/low temperature respectively.  For
the high temperature case, a first rigorous proof was given in
\cite{GalMS}, while our approach is based on the polymer system
introduced in \cite{GrKu}.  For an overview of these methods and references 
we refer to \cite{Simon}.

Clearly the choice of the first approximation $\BARH^{(0)}_M$ is 
crucial to our
method and it should be such that (i) it is explicitly computable as
the constructions in \cite{steph,KMVPNAS}, and (ii) it provides a
good estimate between the microscopic and the 2nd-order 
coarse-graining.  In a subsequent paper \cite{KPRT2} we will
show how to apply these ideas and techniques to systems with both
short and long-range interactions. In this case the choice of the
first approximation $\BARH^{(0)}_M$ will not be given in a closed
form but it will be computed numerically in an efficient way.

All error estimates are calculated in terms of the specific relative
entropy of the corresponding equilibrium Gibbs measures. Specific
relative entropy represents the loss of information in the transition
from the microscopic to the coarse-grained models and here is used to
assess the information compression for the same level of coarse
graining in schemes that are of a higher or lower order, determined by
the truncation level $p$ in (\ref{intro_expansion}).  The primary
practical purpose of the higher-order corrections is in allowing us to
extend the regime of validity of the expansion and to obtain very
accurate coarse grainings of the Gibbs measure even if the parameter
$\epsilon$ in \VIZ{intro_expansion} is not necessarily much smaller 
than one.  We refer to 
examples in Section~\ref{sims} where even for $\epsilon=\BIGO(1)$, we
obtain an accurate prediction of the hysteresis when we include the
higher-order correction terms.  

The error analysis developed in this paper also provides an error
expansion that can be computed and tracked in the process of a
simulation from the numerical data.  Such an {\it a posteriori} error
cannot be numerically computed directly from the relative entropy
formula, since it involves the calculation of the probability density
functions of the microscopic and the coarse-grained measure.  However,
the expansion in (\ref{intro_expansion}) shows a constructive way of
calculating the error made, for instance, by the 2nd-order
coarse-graining $\BARH_M^{(0)}(\eta)$ in terms of other coarse
observables given by the higher-order correction terms $\BARH_M^{(p)}$
plus a controlled error of order $\BIGO(\epsilon^{p+1})$.  
Similarly to the numerical analysis of
approximations for PDEs, we derive a
priori and a posteriori error estimates between the exact microscopic
solution and the approximating coarse-grained one.  In contrast to
the PDE setting where such error analysis is calculated in a suitable
norm, here the error is measured in terms of the specific relative
entropy.
As in the case of PDEs such a posteriori errors are also useful for
constructing adaptive methods. For some earlier
work on adaptivity for stochastic systems see, e.g., \cite{CKV1, CKV2, Szepessy1}.

The paper is organized as follows: In Section~\ref{sec2} we present the main
results and an outline of our methods. In Section~\ref{sec-cluster} we discuss cluster
expansions in the context of coarse-graining and derive the effective
coarse-grained Hamiltonian as an abstract expansion. In Section~\ref{schemes} we
calculate the first terms in the expansion of the Hamiltonian and
formulate two concrete numerical schemes; Scheme~\ref{scheme1} is
second-order accurate, while Scheme~\ref{scheme2} is third-order
accurate. In Section~\ref{sims} we present simulations with these schemes in a
demanding phase transition regime. Appendix A includes the details of
our analytical calculations in Section~\ref{schemes}. Appendix B contains a brief
description of computational  background for the Monte Carlo algorithms used in the
sampling of microscopic and coarse-grained Gibbs states carried out in
Section~\ref{sims}.

%
%
\medskip\noindent{\bf Acknowledgments:}
The research of M.K. was partially supported by DE-FG02-05ER25702,
NSF-DMS-0413864 and
NSF-ITR-0219211. The research of P.P. was partially supported by 
NSF-DMS-0303565.  The research or L.R-B. was partially supported by
NSF-DMS-0306540. The research of D.T. was partially supported by 
NSF-DMS-0413864 and
NSF-ITR-0219211. 

%
%

\section{Main results and outline of the method}\label{sec2}

%
%
\subsection{Microscopic models}
We consider as the physical domain for the system the
$d$-dimen\-sional torus $\T_d:=[0,1)^d$ with periodic boundary
conditions. The microscopic system consists of a uniform lattice
$\LATT:=(\frac{1}{n}\Z)^d\cap\T_d$.  The number of lattice sites
$N=n^d$ is fixed, but arbitrary and finite. We consider here periodic
boundary conditions, but other boundary conditions can be accommodated
easily.

The spin (or order parameter) $\sigma(x)$ takes values in $\{+1,-1\}$
at each lattice site $x\in\LATT$. A spin configuration
$\sigma=\{\sigma(x)\}_{x\in\LATT}$ is an element of the configuration
space $\SIGMA:=\{+1,-1\}^{\LATT}$. The energy of the configuration
$\sigma$ is given by the Hamiltonian
\begin{equation}\label{microHamiltonian}
H_N(\sigma)=-\frac{1}{2}\sum_{x\in\LATT}\sum_{y\neq x}J(x-y)\sigma(x)\sigma(y)+
              \sum_{x\in\LATT}h(x)\sigma(x) \COMMA
\end{equation}
where the two-body inter-particle potential $J$ describes the interaction between individual spins 
and $h$ is an external field.  

The strength of the potential is measured by
$\NORM{J}\equiv\sum_{x\neq 0} |J(x)|$, i.e. we assume that the
two-body potential is summable.  As we will scale later the
Hamiltonian with the inverse temperature $\beta$ we can assume,
without loss of generality, that $\NORM{J}=1$.

\smallskip
\noindent{\em Example 1: Nearest-neighbor interaction.} In this case
the spin at site $x$ interacts only with its nearest neighbors
on the lattice $\LATT$, i.e.
$$
J(x-y)=\begin{cases}
                 J  & \text{if $|x-y| = \frac{1}{n}$;} \\  
                 0 & \text{otherwise.}
       \end{cases}
$$
  
\noindent{\em Example 2: Finite-range interactions.}  A spin at site
$x$ interacts with its neighbors which are at most $\RANGE$ lattice
points away from $x$.  It will be useful to consider the range of the
interaction $\RANGE$ as a parameter of the model. To do this let
\begin{equation}
  \JNOT:\R^+\to\R\COMMA\SPACE 
  \JNOT(r) = 0 \COMMA\SPACE\mbox{ if $|r|\geq 1$.} \label{defJV2} 
\end{equation}
Then the potential $J(x-y)$ can be taken to have the form 
\begin{equation}
  J(x-y) = \frac{1}{\RANGE^d}\JNOT\left(\frac{n}{\RANGE}|x-y|\right)\COMMA
            \SPACE x,y\in\LATT\PERIOD \label{defJV1}\\
\end{equation}
The factor $1/\RANGE^d$ in \eqref{defJV1} is a normalization which
ensures that the strength of the potential $J$ is essentially
independent of $L$ and we have $\NORM{J} \simeq \int |V(r)| dr$.  Note
also that Example 1 can be obtained from Example 2 by setting $L=1$.

\smallskip
\noindent{\em Example 3: Long-range interactions.} In this case we
assume that a spin interacts with all spins on the lattice $\LATT$ via
a summable interaction $J(x)$.  Since our goal is to construct
coarse-grained approximations of the Hamiltonian suitable for
numerical simulations it will be convenient to truncate long-range
interactions and control the error term. In order to do this we choose
a small parameter $\delta$ and choose $\RANGE$ such that
$$
  \sum_{\{x \,:\, |x| \ge \frac{\RANGE}{n}\}} |J(x)| \le \delta \PERIOD    
$$ 
The parameter $\RANGE$ can be thought as the effective range of the
potential and we can truncate the long-range potential $J(x)$ by
setting $J(x)=0$ if $|x| \ge \frac{\RANGE}{n}$.  If we denote by
$\hat{J}$ the truncated potential and by $\hat{H}_{N}$ the
corresponding Hamiltonian it is easy to see that
$$
  \frac{1}{N} | H_N(\sigma) - \hat{H}_N(\sigma)| \,=\,  O(\delta) \PERIOD
$$
i.e., the error per unit volume is of order $\delta$.  
We can then assume that the truncated potential 
has the form \eqref{defJV1}. 

\smallskip
\noindent{\em Example 4: Kac-type interactions.}  Mean-field
interactions are obtained formally by taking a finite range
interaction and setting $L=n$ in \eqref{defJV1}
\begin{equation}
   J(x-y) = \frac{1}{N}\JNOT\left(|x-y|\right)\COMMA
            \SPACE x,y\in\LATT\PERIOD \label{defJV3}\\
\end{equation}
Note that this is different from a summable long-range interaction since
the potential is scaled with the size of the system.

\smallskip
The potential in Example 2 will be central to our analysis.  Both
Example 1 and 4 are obtained in suitable limiting cases ($\RANGE=1$
and $\RANGE=n$ respectively) and the potential in Example 3 can be
approximated by such potentials (see the error estimates below).

The finite-volume equilibrium states of the system are given by the
canonical Gibbs measure
\begin{equation}\label{microGibbs}
    \mu_{N,\beta}(d\sigma)=\frac{1}{Z_N}\EXP{-\beta H_N(\sigma)}
                            P_N(d\sigma)\COMMA
\end{equation}
where $\beta$ is the inverse temperature, $Z_N$ is the normalizing
partition function, and $P_N(d\sigma)$, the prior distribution on
$\SIGMA$, is the product measure
$$
   P_N(d\sigma)=\prod_{x\in\LATT}\rho(d\sigma(x))\PERIOD
$$
A typical choice $\rho(\sigma(x)=+1)=\frac{1}{2}$ and $\rho(\sigma(x)=-1)=\frac{1}{2}$ 
is the distribution of a Bernoulli random variable for each $x\in\LATT$.
%
%
%
\subsection{Relative specific entropy} 

One of the obvious issues arising in any attempt to coarse-grain
microscopic systems is the evaluation of the numerical error as we
move from finer to coarser scales. This error essentially involves the
comparison of the microscopic and the coarse-grained probability
measures.  A natural way to compare two probability measures is to
compute their relative entropy.  Let $\pi_1(\sigma)$ and
$\pi_2(\sigma)$ be two probability measures defined on a common
probability space, with the state space $\STATESP$.  The relative
entropy, or Kullback-Leibler distance, of $\pi_1$ with respect to
$\pi_2$, is defined as
\begin{equation}\label{relent}
  \RELENT{\pi_1}{\pi_2} \equiv \int_{\STATESP} 
        \log\left( \frac{d\pi_1}{d\pi_2}\right) d\pi_1 \COMMA
\end{equation}
where $d\pi_1/d\pi_2$ is the Radon-Nikodym derivative of $\pi_1$ with
respect to $\pi_2$, On a countable state space $\STATESP$ we obtain
$$
  \RELENT{\pi_1}{\pi_2} = \sum_{\sigma\in\STATESP} \pi_1(\sigma) 
                          \log\frac{\pi_1(\sigma)}{\pi_2(\sigma)}\PERIOD
$$ In information theory the relative entropy $\RELENT{\pi_1}{\pi_2}$
provides a measure of ``information loss'' or ``information distance''
of $\pi_1$ compared to $\pi_2$. In our context we use the relative
entropy in order to assess the information compression of different
coarse-graining schemes. Basic properties and applications in
information theory can be found in \cite{CoverThomas}.

Using Jensen's inequality it is not difficult to show that
\begin{eqnarray*}
  \RELENT{\pi_1}{\pi_2} &\geq& 0\;\;\;\mbox{and,} \\
  \RELENT{\pi_1}{\pi_2} & = & 0\;\;\;
           \mbox{ if and only if $\pi_1(\sigma) = \pi_2(\sigma)$ for all $\sigma\in \STATESP$.} \\
  \RELENT{\pi_1}{\pi_2} & = & \infty \;\;\; \mbox{if $\pi_1$ is not absolutely continuous with
                                                  respect to $\pi_2$.} 
\end{eqnarray*}
The relative entropy is not a metric (it is not symmetric), but one
can use the relative entropy to bound the total variation distance of
the measures $\pi_1$ and $\pi_2$ as demonstrated by 
{\it Csisz\'ar-Kullback-Pinsker inequality}
(\cite{CoverThomas})
\begin{equation} \label{totvar}
  \RELENT{\pi_1}{\pi_2} \geq \frac{1}{2} \left(\sum_{\sigma\in \STATESP} 
               |\pi_1(\sigma) - \pi_2(\sigma)| \right)^2 \equiv 
               \frac{1}{2} \|\pi_1 - \pi_2\|^2_{\mathrm{TV}}\COMMA
\end{equation}
or equivalently using the dual form of the total variation norm
\begin{equation}
   \sup_{\|f\|_\infty =1} \left| \int f \,d\pi_1 - \int f \,d\pi_2\right| 
   \leq \sqrt{\RELENT{\pi_1}{\pi_2}}\COMMA
\end{equation} 
where $\| \cdot \|_\infty$ denotes the usual sup-norm. 

In statistical mechanics, the entropy, which is the relative entropy
with respect to the prior distribution, and the relative entropy play
a prominent role, in particular in the variational principles. The
relative entropy is, in general, difficult to compute since it
requires the sampling of the probability distribution $\pi_1(\sigma)$
and $\pi_2(\sigma)$ which can be prohibitively expensive if the
dimension of the state space is large. For Gibbs measures it is,
however, a little simpler in the sense that it depends only on the
partition functions and the expected values of the Hamiltonian.
Consider, for example, two Gibbs measures $\mu^{(1)}_{N,\beta}$ and
$\mu^{(2)}_{N,\beta}$ with Hamiltonians $H^{(1)}_N$ and $H^{(2)}_N$
and partition functions $Z^{(1)}_N$ and $Z^{(2)}_N$. Then we have
\begin{equation} \label{relentgibbs}
\RELENT{\mu^{(1)}_{N,\beta}}{\mu^{(2)}_{N,\beta}} \,=\,\log \left( 
\frac{Z^{(2)}_N}{Z^{(1)}_N}\right) + \int 
\left( H^{(2)}_N - H^{(1)}_N \right) d\mu^{(1)}_{N,\beta} \PERIOD   
\end{equation}
Note that the Hamiltonians and the logarithm of partition functions
are extensive quantities: they are proportional to the size of the
system $N$.  The formula \eqref{relentgibbs} shows that the entropy of
a Gibbs measure and the relative entropy of a Gibbs measure with
respect to another Gibbs measure are also extensive quantities.  It is
therefore natural to compare two Gibbs measures on $\SIGMA$ by
computing their specific relative entropy
$$
  \frac{1}{N} \RELENT{\mu^{(1)}_{N,\beta}}{\mu^{(2)}_{N,\beta}} \,.
$$ The specific relative entropy is a measure of the information loss
per unit volume and will be our main estimation tool in this paper.  A
simple example were this interpretation is evident arises if one
considers $N$ independent particles and the ensuing scaling of the
collective discrepancy if an individual error is committed on the
observation of each particle.
%
%
%
%
\subsection{Coarse-grained models}  
The coarse-graining procedure consists of three steps that we describe
separately.

\noindent 
(a) {\em Coarse graining of the configuration space.}  We partition the torus 
$\T_d$ into $M=m^d$ cells: For $k=(k_1,\ldots,k_d)\in\Z^d$ with $0 \le k_i 
\le m-1$ we define
$\CUBE_k\equiv [\frac{k_1}{m},\frac{k_1+1}{m})\times\cdots\times
[\frac{k_d}{m},\frac{k_d+1}{m})$ and we have $\T_d = \cup_{k}
\CUBE_k$.  We identify each cell $\CUBE_k$ with a lattice point of the
coarse lattice $\LATTC = (\frac{1}{m}\Z)^d\cap\T_d$.  Each coarse cell
contains $Q=q^d$ points of the microscopic lattice points with
$N\equiv n^d=(m q)^d\equiv MQ$.  We will refer to $Q$ as the level of
coarse graining ($Q=1$ corresponds to no coarse graining).

We assign a new spin value $\eta(k)$ for the cell $\CUBE_k$ according
to the rule
$$
   \eta(k)=\sum_{x\in \CUBE_k}\sigma(x)\PERIOD
$$ The spin $\eta(k)$ takes values in $\{ -Q, -Q+2, \ldots, Q\}$ and the
configuration space for the coarse grained system is $\SIGMAC\equiv \{
-Q, -Q+2, \ldots, Q\}^{\LATTC}$.  We denote by $\COP$ the map
$$
  \COP :\SIGMA \rightarrow \SIGMAC\COMMA\;\;\;\;
       \sigma \mapsto \{\sum_{x\in \CUBE_k}\sigma(x)\}_k \COMMA
$$ which assigns a configuration $\eta = \{\eta(k)\}_{k\in\LATTC}$ on
the coarse lattice given a configuration $\sigma \,=\,
\{\sigma(x)\}_{x\in\LATT}$.  An equivalent coarse grained variable,
which we shall also use later, is
$$ 
  \alfa(k):=\card \{ x\in\CUBE_k:\,\sigma(x)=+1\}
$$ 
which takes values in $\{0,1,\ldots, \CGPARQ\}$.  
Both coarse variables are equivalent 
and they are related by the transformation 
$\eta=2\alfa-\CGPARQ$  or $\alfa=\frac{\eta+\CGPARQ}{2}$.

\smallskip
\noindent
(b) {\em Coarse-graining of the prior distribution.} The prior distribution $P_N$ on 
$\SIGMA$ induces a new prior distribution on $ \SIGMAC$ given by ${\bar P}_M = 
P_N \circ \COP^{-1}$, i.e., 
$$
  \BARIT{P}_M(\eta) = P_N(\sigma \,:\, \COP(\sigma) = \eta)\PERIOD
$$ 
Since $\eta(k)$ depends only on the spin $\sigma(x)$, with
$x \in \CUBE_k$, the measure ${\bar P}_M$ is a product measure
$$ 
  \BARIT{P}_M( d{\eta})=  \prod_{k \in \LATTC} {\bar \rho}(d \eta(k)) \COMMA
$$
with  
$$ 
  \BARIT{\rho}(\eta(k) ) \,=\, \binom{Q}{\frac{\eta(k)+ Q}{2}}\left(\frac{1}{2}\right)^{Q}
  \PERIOD
$$

The conditional probability  $P_N(d\sigma|{\eta})$ plays a crucial role in the sequel.  
Since $\eta(k)$ depends only on the spin $\sigma(x)$ with $x \in \CUBE_k$,  the probability 
$P_N(d\sigma|{\eta})$ factorizes over the coarse cells. We have
\begin{equation}\label{cond_meas}
   P_N(d\sigma|{\eta})=\frac{ P_N( \sigma\cap\{\COP(\sigma)=\eta\})}{\BARIT{P}_M(\eta)}
                  =\prod_{k \in \LATTC} \tilde{\rho}_{k,\eta(k)}(d\sigma) \,.
\end{equation}
and $\tilde{\rho}_{k,\eta(k)}(d\sigma)$ depends only on $\{\sigma(x)\}_{x \in \CUBE_k}$.   
In particular we have 
\begin{equation}
  \tilde{\rho}_{k,\eta(k)}(\sigma(x)=1)=\frac{\eta(k) + Q}{2 Q}\,\,\mbox{and}\,\,\,
  \tilde{\rho}_{k,\eta(k)}(\sigma(x)=-1)=\frac{Q-\eta(k)}{2 Q}\PERIOD
\end{equation}
To simplify the notation and because our estimates are uniform in $\eta(k)$ 
we denote this measure simply by $\tilde{\rho}_k$. 
For a function $f=f(\sigma)$ we define the conditional  expectation 
\begin{eqnarray}
  \EXPECT[  f |{\eta}] \,&=&\, \int f(\sigma) P_N(d\sigma|{\eta}) \,=\, 
  \frac{1}{\BARIT{P}_M(\eta)} \int_{\{\sigma:\, \COP(\sigma)={\eta}\}} f(\sigma) P_N(d\sigma)  
\nonumber \\
\,&=&\,  \int f(\sigma) \prod_k \tilde{\rho}_k (d\sigma) \,.
\end{eqnarray}

\smallskip
\noindent
(c) {\em Coarse-graining of the Hamiltonian.}  We want to construct a Hamiltonian 
$\BARH_M(\VEC{\eta})$ at the coarse-level.  A natural definition of such Hamiltonian, as we explain below, is 
given by the renormalization group block averaging transformation (also known as  
Kadanoff transformation). 
\begin{definition}
   The {\em exact} coarse grained Hamiltonian 
   $\BARH_M(\VEC{\eta})$ is defined by the formula
   \begin{eqnarray}\label{rg}
     \EXP{-\beta \BARH_M (\eta)} & \,= \,&\EXPECT[\EXP{-\beta H_N}|{\eta}]\PERIOD
   \end{eqnarray}
   Given the Hamiltonian $\BARH_M$ we define 
   the corresponding Gibbs measure by
   \begin{equation}\label{cg_gibbs}
     \BARIT{\mu}_{M,\beta}(d{\eta})=\frac{1}{\BARIT{Z}_M}
     \EXP{-\beta\BARH_M({\eta})} \BARIT{P}_M( d{\eta})\PERIOD
   \end{equation}
\end{definition}
The factor $\beta$ in front of $\BARH_M(\eta)$ is merely a convention, 
in general the Hamiltonian $\BARH_M(\eta)$ does itself depend on $\beta$ in a nonlinear, complicated way.  

Clearly for a fixed value of $N$, this transformation is well defined. 
However, from the point of view of statistical mechanics this is not sufficient. 
We want to construct the coarse-grained Hamiltonian
$\BARH_M(\eta)$, for any $M$,  as a sum of (summable) many-body interactions.  
Our original Hamiltonian $H_N(\sigma)$ involves only one-body and two-body interactions but the 
nonlinear transformation \eqref{rg} will not preserve this property.   
The Hamiltonian we will consider in the sequel will have the form 
\begin{equation} 
K_M(\eta) \,=\, \sum_{X\subset \LATTC} J(X, \{\eta_k\}_{k \in X}) 
\end{equation}
where the $J(X, \{\eta_k\}_{k \in X})$ are translation-invariant many-body interactions 
involving ${\rm card(X)}$-many different sites.  Such an interaction is said to be summable if 
$$
\sum_{\{X\,:\,0\in X\}} \|J(X, \{\eta_k\}_{k \in X})\|<\infty \PERIOD
$$ 
 
It is well-known \cite{EFS} that the Kadanoff transformation suffers from some 
pathologies. At very low temperature the Kadanoff transformation for the nearest-neighbor 
Ising model with the zero magnetic field is not well-defined in the sense that the coarse-grained 
Hamiltonian $\BARH_M$ does not correspond to a summable interaction.  These pathologies are  
relatively mild \cite{BKL} and can be eliminated by extending slightly the concept of Gibbs measures. 
For a comprehensive presentation of this issue see the review paper \cite{EFS}
or also \cite{BKL} among others.
Furthermore, these pathologies will not play any role in our 
analysis: the cluster expansion techniques we are using exclude the occurrence of these 
pathologies for the models and the values of the parameter we consider. 

The exactness of the coarse graining in the Kadanoff transform
is expressed by the fact that the
specific relative entropy of $\BARIT{\mu}_{M,\beta}(d{\eta})$ with respect 
to the coarse-grained Gibbs measure $\mu_{N,\beta} \circ\COP^{-1}$ vanishes. 
Indeed we have
\begin{eqnarray}
   Z_N  &=& \int \EXP{-\beta H_N(\sigma)}P_N(d\sigma) = \int 
            \EXPECT[\EXP{-\beta H_N}|{\eta}] \, \BARIT{P}_M(d{\eta}) \nonumber \\
        &=& \int \EXP{-\beta \BARH_M({\eta})}  \, \BARIT{P}_M(d{\eta}) = \BARIT{Z}_M\COMMA
\end{eqnarray}
and consequently from \VIZ{rg},
\[
   \frac{1}{\BARIT{Z}_M}\EXP{-\beta\BARH_M(\eta)} \BARIT{P}_M(\eta)
  =\frac{1}{Z_N}\int_{\{\sigma\SEP\, \COP(\sigma)={\eta}\}}\EXP{-\beta H_N} P_N(d\sigma)\COMMA
\]
and thus
$$
   \frac{1}{N}\RELENTR(\BARIT{\mu}_{M,\beta}|\mu_{N,\beta}\circ\COP^{-1}) \,=\, 0 \PERIOD
$$ 

Even for moderately large values of $N$ the exact computation of
$\BARH_M(\eta)$ is, in general, impractical. The coarse grained
Hamiltonian involves not only two-body interactions, but also many-body
interactions of an arbitrary number of spins.  Our goal is to present a systematic way 
of calculating explicit approximations of the coarse-grained Hamiltonian $\BARH_M$, to any
given degree of accuracy.
In the first step we define an approximate
Hamiltonian $\BARH^{(0)}_M(\eta)$ and we give an a priori
bound on the blocking error. 
The choice of $\BARH^{(0)}_M$ we make, following \cite{KMV}, is  given in
the following definition.
\begin{definition}\label{hbar0}
  The first approximation $\BARH^{(0)}_M$ of the coarse-grained Hamiltonian
  $\BARH_M$ is given by the formula
  \begin{equation} \label{cg_Hamiltonian}
     \BARH^{(0)}_M(\eta)\,\equiv\, \EXPECT[H_N|\eta] \PERIOD
  \end{equation}
\end{definition}
In order to compute $\BARH^{(0)}_M(\eta)$ we note that we have
\begin{eqnarray*}
&&  \EXPECT[ \sigma(x)|\eta]=  \frac{\eta(k)}{Q}\COMMA  \SPACE  x \in \CUBE_k  \COMMA \\
&&  \EXPECT[ \sigma(x)\sigma(y)|\eta] = \frac{\eta(k)^2 - Q}{Q(Q-1)}\COMMA  
             \SPACE  x \in \CUBE_k,  y \in \CUBE_k  \PERIOD
\end{eqnarray*}
Then
$$ 
  \BARH^{(0)}_M(\eta) \,=\, -\frac{1}{2}\sum_{k}\sum_{l\neq k}\bar{J}(k-l)\eta(k)\eta(l)
                            -\frac{1}{2}\sum_{k}\bar{J}(0)(\eta^2(k)-Q)+ h\sum_k \eta(k) 
$$
where
\begin{eqnarray*}
&&  \BARJ(k-l)\,=\, \frac{1}{Q^{2}}\sum_{x\in \CUBE_k,y\in \CUBE_l}J(x-y)\COMMA
    \;\;\;\mbox{for $k\neq l$,}\\
&& \BARJ(0)=\frac{1}{Q(Q-1)}\sum_{x,y\in \CUBE_k,y\neq x}J(x-y)\COMMA
    \;\;\;\mbox{for $k=l$.}
\end{eqnarray*}


By defining $\BARH^{(0)}_M(\eta)$  one 
replaces the potential $J(x-y)$ by its average over a coarse cell. 
Thus the error for the potential is proportional to 
$$
E_{kl}(x-y):=J(x-y)-\bar{J}(k,l)\COMMA \SPACE x \in \CUBE_k, y \in \CUBE_l\COMMA
$$
which measures the variation of the potential $J(x-y)$ over a cell.  
An estimate on the error is provided by the following lemma

\begin{lemma}\label{a-priori_error1}
   Assume that $J$ satisfies \eqref{defJV2}--\eqref{defJV1} and assume that $V(r)$ is ${\cal C}^1$. 
  \begin{enumerate}
     \item There exists a constant $C>0$ such that, if $x\in \CUBE_k$ and $y\in\CUBE_l$,  we have
          \begin{equation}\label{errJb1}
		|J(x-y) - \bar J(k,l)| \leq 2\frac{q}{L^{d+1}}
                \sup_{\genfrac{}{}{0cm}{2}{x'\in \CUBE_k,}{y'\in \CUBE_l}}\|\nabla V(x'-y')\| \PERIOD
          \end{equation}
     \item There exists a constant $C>0$ such that, if $\COP(\sigma)=\eta$, we have 
          \begin{equation}
                \frac{1}{N} \left| H_N(\sigma) - \BARH_M^{(0)}(\eta)\right| \,\le\, 
                        C \frac{q}{L} \|\nabla V\|_\infty \PERIOD
          \end{equation}
  \end{enumerate}
\end{lemma}
The lemma will be proved in Subsection~\ref{errors}.  
Note that the estimates in Lemma~\ref{a-priori_error1} 
are not necessarily optimal and can be improved under suitable assumptions on the potential.  
The importance of Lemma~\ref{a-priori_error1} lies in the identification of the 
small parameter in Theorem~\ref{mainresult} below 
\begin{equation}\label{small_parameter}
   \epsilon \,\equiv\,   C \beta \frac{q}{L} \|\nabla V\|_\infty \COMMA
\end{equation}
where $C$ is some constant.
Note that we included the inverse temperature $\beta$ in the parameter, taking into account that  
$\beta$ multiplies the Hamiltonian in the Gibbs measure.  The parameter $\epsilon$ 
encapsulates the various factors involved in our coarse-graining method. 

\noindent{\it (i) Coarse-graining cell size.} The factor $q/L$ in \VIZ{small_parameter} indicates 
how the size of the coarse cell governs the effectiveness of the coarse-graining procedure. 

\noindent{\it (ii) Slow variation of the potential.} The factor $\|\nabla V\|_\infty$ in \ref{small_parameter} 
indicates that the coarse-graining will be particularly effective even for large cell sizes if the 
potential $V$ is slowly varying. This fact is exemplified by considering the extreme case of the
Curie-Weiss model where $J(x-y)=J$ is  independent of $x-y$.  It is easy to 
see that in this case the coarse-graining is actually {\em exact}. 

\noindent{\it (iii) Temperature.} At high temperature, i.e, for small
$\beta$ the coarse-graining procedure will work well even for large
values of $q/L$, i.e., for very large coarse cells or for potentials
with very short range.  This reflects the fact that at high
temperature Gibbs measures are actually very close to product measures
and that, trivially, product measures can be coarse-grained exactly.
On the other hand our methods clearly do not apply to very-low
temperatures.

Let us denote by $\BARIT{\mu}_{M,\beta}^{(0)}(d{\eta})$ the Gibbs
measure
\begin{equation*}
   \BARIT{\mu}_{M,\beta}^{(0)}(d{\eta}) \,=\, \frac{1}{\BARIT{Z}_M^{(0)}} \EXP{ -\beta  \BARH^{(0)}_M(\eta)}  
   \BARIT{P}_M (d{\eta}) \PERIOD
\end{equation*}
From Lemma~\ref{a-priori_error1} we obtain immediately the bound
\begin{equation}\label{errb00}
 \frac{1}{N}\RELENTR(\BARIT{\mu}^{(0)}_{M,\beta}|\mu_{N,\beta}\circ\COP^{-1}) \,=\,  O(\epsilon) \,.
\end{equation} 
Furthermore, Lemma~\ref{a-priori_error} allows us to perform a cluster expansion and compute 
higher-order corrections.  The basic idea is to rewrite the exact coarse-graining as
$$
  \EXP{-\beta\BARH_M(\eta)}=\EXP{-\beta\BARH_M^{(0)}(\eta)}\EXPECT[\EXP{-\beta
  (H_N(\sigma)-\BARH^{(0)}_M(\eta))}|\eta] \COMMA
$$
or 
\begin{equation}\label{perturbation}
  \BARH_M(\eta)=\BARH^{(0)}_M(\eta)-\frac{1}{\beta}
                \log\EXPECT[\EXP{-\beta(H_N(\sigma)-\BARH^{(0)}_M(\eta))}|\eta]\PERIOD
\end{equation}
Note that the exponential in the second term is not necessarily small, it is 
in fact of order $N\epsilon$.  
Cluster  expansions are tools which allows to expand such quantities in convergent power series 
using the independence properties of product measures.  The crucial fact here is that 
conditional measure  $P_N(d\sigma|{\eta})$ factorizes over the coarse cells.  
The main result of our paper, proved in Section~\ref{schemes}, is 
\begin{theorem}\label{mainresult} 
Let us assume that $J$ satisfies \eqref{defJV2}--\eqref{defJV1} and that 
$V(r)$ is ${\cal C}^1$.   Then there exists a constant $\delta_0>0$ such that if 
$$
  \delta= Q \epsilon < \delta_0 \COMMA
$$
the  Hamiltonian $\BARH_M(\eta)$ can be expanded into a convergent series
$$
 \BARH_M(\eta)=\sum_{p=0}^{\infty} \BARH^{(p)}_M(\eta) \COMMA
$$
where each term $\BARH^{(p)}_M(\eta)$ is a sum of finite-range translation invariant many-body potentials. 
The first few terms are explicitly calculated in Scheme~\ref{scheme1} and Scheme~\ref{scheme2} 
in Section~\ref{schemes}.
We have the following error bounds uniformly in $\eta$ and $N$
$$ 
   \frac{\beta}{N} \left( \BARH_M(\eta)-(\BARH^{(0)}_M(\eta) + \ldots +
   \BARH^{(p)}_M(\eta) \right) \,=\, \BIGO(\epsilon^{p+1})\PERIOD
$$
If we define the Gibbs measures  
\begin{equation*}
   \BARIT{\mu}_{M,\beta}^{(p)}(d{\eta}) \,=\, \frac{1}{\BARIT{Z}_M^{(p)}} \EXP{ -\beta( \BARH^{(0)}_M(\eta) 
      + \ldots + \BARH^{(p)}_M(\eta))}  \BARIT{P}_M (d{\eta}) \PERIOD
\end{equation*}
then the following bounds for the relative entropy per unit volume hold
\begin{eqnarray*}
&& \frac{1}{N}\RELENTR(\BARIT{\mu}^{(0)}_{M,\beta}|\mu_{N,\beta}\circ\COP^{-1}) 
    \,=\, \BIGO(\epsilon^2) \COMMA \label{errb0} \\
&& \frac{1}{N} \RELENTR(\BARIT{\mu}_{M,\beta}^{(p)}|\mu_{N,\beta}\circ \COP^{-1}) = 
    \BIGO(\epsilon^{p+1})\COMMA \label{errb}
\end{eqnarray*}
where $p=2,\ldots$ and $\epsilon$ is given by (\ref{small_parameter}).
\end{theorem}
  
We note that in  information theory the relative entropy  provides 
a measure of  ``information distance''  of two probability measures. 
In our context we  use the relative entropy estimate of Theorem~\ref{mainresult} in order to assess
the information compression of different coarse-graining schemes.
According to how many correction terms are included in the expansion corresponding 
to different truncation levels $p$, such as Scheme~\ref{scheme1} and Scheme~\ref{scheme2} 
in Section~\ref{schemes}, we quantify the amount of the information loss when the measures
are compared at the same level of coarse graining.

\begin{remark}
{\rm
\begin{enumerate}
\item The cluster expansion  provides an explicit  algorithm to 
      compute the higher-order approximations $\BARH^{(p)}_M$.   
      For example, the next order correction is given by $\BARH^{(1)}_M+\BARH^{(2)}_M$ as in
      Scheme~\ref{scheme2} in Section~\ref{schemes}. The analytical computations required for  determining 
      higher-order corrections become quickly very involved, nonetheless they can be carried out
      with symbolic computational tools. 
\item The primary practical purpose of the  higher-order estimates
      in Theorem~\ref{mainresult} is  in allowing us to extend the regime of
      validity of the expansion and obtain very accurate coarse-graining of the Gibbs measure even if
      the parameter $\epsilon$ given by  (\ref{small_parameter}) is not necessarily much 
      smaller than one. We refer to
      the examples in Section~\ref{sims} where even for $q\geq L$, and $\beta>1$, i.e., $\epsilon=\BIGO(1)$, 
      we get an accurate prediction of critical behavior  when
      we include the higher-order correction terms. 
\item The reader might wonder about the appearance of a new small parameter $\delta=Q\epsilon$ 
      in the statement of Theorem~\ref{mainresult}.  The reason lies in the details of 
      the cluster expansion. Since the conditional measure $P_N(d\sigma|{\eta})$ factorizes 
      over the coarse cells, it is required for the cluster expansion to converge that the error 
      in every coarse cell is small. The error for each site of the microscopic lattice is of 
      order $\epsilon$ and thus the error 
      for a coarse cell is of order $\delta=Q\epsilon$.  
\item Note further that the bound (\ref{errb0}) improves on the bound (\ref{errb00}) and that
      shows the first approximation $\BARH_M^{(0)}$ is actually already a second-order method. 
      This is due to cancellations 
      which follow directly from the definition of $\BARH_M^{(0)}$, see Definition~\ref{hbar0}
      and \VIZ{deltaJiszero} below.
\end{enumerate}
}
\end{remark}

In the previous theorem the errors are calculated with respect to the
relative entropy of the corresponding equilibrium measures.
Apart from this a priori estimate, we would also like to 
have a formulation of the error that can be explicitly computed,
from the data of the simulation.
In our case, such an {\it a posteriori} error cannot be numerically computed
directly from the relative entropy formula, 
since it involves the calculation of the probability density functions
of the microscopic and the coarse-grained measure.
However,  the error estimate of 
Theorem~\ref{mainresult} provides us with an explicit way of calculating the error
made, for instance, by the 2nd-order coarse-graining.
In Section~\ref{schemes} we prove the a posteriori error estimate
for $\RELENTR(\BARIT{\mu}_{M,\beta}^{(0)}|\mu_{N,\beta}\circ \COP^{-1})$
given in the theorem
\begin{theorem}
[{\rm A posteriori error}]
\label{mainresult2}
We have
\[
  \RELENTR(\BARIT{\mu}_{M,\beta}^{(0)}|\mu_{N,\beta}\circ \COP^{-1})=
  \EXPECT_{\BARIT{\mu}_{M,\beta}^{(0)}}[R(\eta)]+
  \log\left(\EXPECT_{\BARIT{\mu}_{M,\beta}^{(0)}}[\EXP{R(\eta)}]\right)
  +\BIGO(\epsilon^3)\COMMA
\]
where the residuum operator $R(.)$ is given by 
$$
R(\eta)=\BARH^{(1)}_M(\eta)+\BARH^{(2)}_M(\eta)\PERIOD
$$
\end{theorem}
Note  that the quantity $\EXPECT_{\BARIT{\mu}_{M,\beta}^{(0)}}[R(\eta)]$
can be calculated on-the-fly using the coarse simulation.
Calculations involving the a posteriori error estimation and  related
adaptive methods will be discussed in a forthcoming publication.
Earlier work that uses only an upper bound and not the sharp estimate 
of Theorem~\ref{mainresult2}
can be found in \cite{CKV1,CKV2}. 

%
%
\section{Cluster expansion and effective interactions}\label{sec-cluster}
In this section we expand the term 
$\EXPECT[\EXP{-\beta(H_N(\sigma)-\BARH_M^{(0)}(\eta))}|\eta]$ in \VIZ{perturbation}
into a convergent 
series using a cluster expansion.   
\subsection{Introduction to the polymer model}
It is convenient to choose an ordering on the $d$-dimensional lattice $\LATTC$.
For example, the lexicographic ordering defines for $d=2$ the relation $(k_1,k_2)\leq(l_1,l_2)$ 
if and only if $k_1<l_1$ or $k_1=l_1$ and $k_2\leq l_2$ and then extended recursively
for arbitrary dimension $d$.
For later use, we also denote by $B_r(k)$ the ball centered at $k$ of radius $r$
for a given metric on the lattice.
Furthermore, by $l\in B^+_r(k)$ we mean that $l\in B_r(k)\cap\{l>k\}$.
To set-up the cluster expansion we write the difference 
$H_N(\sigma)-\BARH_M^{(0)}(\eta)$ as
\begin{eqnarray}
&& H_N(\sigma)-\BARH_M^{(0)}(\eta)= \sum_{k\leq l} \Delta_{kl}J(\sigma)\COMMA
   \;\;\;\;\mbox{where} \\ \nonumber
&& \Delta_{kl}J(\sigma):= -\HALF\sum_{\genfrac{}{}{0cm}{2}{x\in \CUBE_k}{y\in \CUBE_l,y\neq x}}
   (J(x-y)-\bar{J}(k,l))\sigma(x)\sigma(y)(2-\delta_{kl})\PERIOD \label{star}
\end{eqnarray}
Note that we have excluded the external field contribution
to the Hamiltonian, the role of which will be presented
as a remark at the end of this section.
In order to expand the exponential $\EXP{-\beta(H_N(\sigma)-\BARH^{(0)}_M(\eta))}$ 
into a series we set 
\begin{equation}\label{f}
f_{kl}^{}(\sigma):=\EXP{-\beta\Delta_{kl}J(\sigma)}-1
\end{equation}
and obtain, using the factorization properties of $P_N(d\sigma|\eta)$ 
\begin{equation}\label{first}
\EXPECT[\EXP{-\beta(H_N(\sigma)-\BARH_M^{(0)}(\eta)}|\eta] =
\int\prod_{k\leq l}(1+f_{kl}) \prod_{k}\tilde{\rho}_k(d\sigma) \PERIOD
\end{equation}

We give a short description of the polymer model that we
will use in order to organize the cluster expansion.
The product $\prod_{k\leq l}(1+f_{kl}(\sigma))$
can be expressed as
$\sum_{G\in\Gg_M}\prod_{\{k,l\}\in G}f_{kl}$,
where $\Gg_M$ is a collection of all simple graphs on
$M=m^d$ vertices, i.e., graphs
where each edge appears only once, with additional loops
$({k,k})$ on the same vertex.
A one-dimensional example of an element $G$ is $G=\{\{1,2\},\{2,2\},\{2,5\},\{3,6\}\}$
which corresponds to the term $f_{12}f_{22}f_{25}f_{36}$ in the above sum.
We observe further that each $G$ can be divided into connected 
non-intersecting components
called polymers, where by non-intersecting we mean that they do not share
the same coarse cells.
We define the support $\supp(\gamma)$ of a polymer $\gamma$ 
the set of the coarse cells participating in it.
In the previous example, $G=(\gamma_1,\gamma_2)$ where
$\gamma_1=\{\{1,2\},\{2,2\},\{2,5\}\}$ and $\gamma_2=\{\{3,6\}\}$.
We also say $\gamma_1$ is incompatible with $\gamma_2$ and we write
$\gamma_1\not\sim\gamma_2$ if $\supp(\gamma_1)\cap\supp(\gamma_2)=\emptyset$.
Then \VIZ{first} is equal to
\begin{multline}
  \int\sum_{n\geq 0}\frac{1}{n!}
  \sum_{\genfrac{}{}{0cm}{2}{(\gamma_1,\ldots,\gamma_n)}{\gamma_i\not\sim\gamma_j,i\neq j}}
  \prod_{i=1}^n\prod_{\{k,l\}\in\gamma_i}f_{kl}
  \prod_k\tilde{\rho}_k(d\sigma)= \\
  \sum_{n\geq 0}\frac{1}{n!}
  \sum_{\genfrac{}{}{0cm}{2}{(\gamma_1,\ldots,\gamma_n)}{\gamma_i\not\sim\gamma_j,i\neq j}}
  \prod_{i=1}^n\int\prod_{\{k,l\}\in\gamma_i}f_{kl}
  \prod_{\{k\}\in\supp(\gamma_i)}\tilde{\rho}_k(d\sigma)\COMMA
\end{multline}
since the polymers $\gamma_1,\ldots,\gamma_n$ do not share the same coarse cells.
Note that the factor $\frac{1}{n!}$ takes into account the fact that by
relabeling the connected polymers $\gamma_1,\ldots,\gamma_n$
we get the same $G$.

An equivalent formulation is to consider as our building
blocks the more fundamental (than the polymers) 
clusters $R\subset\LATTC$ (i.e. the set of vertices appearing in each
polymer). For example, the corresponding cluster to $\gamma_1$ is $R_1=\{1,2,5\}$.
Note that to each cluster there are several polymers that correspond to it.
To prevent confusion we warn the reader 
that what we call clusters $R$,
it is called polymers  in \cite{Cammarota}. 
We use the name cluster to distinguish it from the previous 
description of vertices with edges, that we have called polymers.
Using clusters the expansion becomes
\begin{eqnarray}
&& \sum_{n\geq 0}\frac{1}{n!}\sum_{
   \genfrac{}{}{0cm}{2}{(R_1,\ldots,R_n)\in\Rr^n}{R_i\cap R_j=\emptyset,i\neq j}}
   \prod_{i=1}^n \zeta(R_i)\COMMA \;\mbox{where}\nonumber \\
&& \zeta(R)=\int\sum_{g\in G_R}\prod_{\{k,l\}\in g}f_{kl}(\sigma)
   \prod_{\{k\}\in R}\tilde{\rho}_k(d\sigma) \label{activity}
\end{eqnarray}
is called the {\it activity} of the cluster $R$.
By $\Rr$ we denote the space of clusters
$R\subset\LATTC$ with $|R|=\card(R)\geq 1$
and $G_R$ is the set of generalized connected graphs on the
set $R$. A generalized connected graph on $R$ is
a collection of distinct subsets of $R$, in which case the 
vertices of the generalized graph are the
indices of the enumeration of the subsets, and
any two subsets of the vertices of the graph are linked
to each other.
The edges of this graph are defined whenever the distinct
subsets share a common element.
For more details we refer to \cite{PdLS}.
Thus we are in the context of \cite[Theorem 2]{Cammarota} and 
we can formulate our first lemma.

\begin{lemma}[{\rm Expansion of the Hamiltonian}]\label{lemma1}
If for every integer $r\geq 2$, $\zeta$ satisfies the condition
\begin{equation}\label{condition}
  \sup_{k\in\LATTC}\sum_{
  \genfrac{}{}{0cm}{2}{R\in\Rr,R\supset \{k\}}{|R|=r}}|\zeta(R)|\leq\delta^{r-1}
\end{equation}
for some $\delta$ such that
$\sup_{k\in\LATTC}\sum_{
\genfrac{}{}{0cm}{2}{R\in\Rr,R\supset \{k\}}{|R|=1}}|\zeta(R)|\leq\delta$
and $\delta<\frac{1}{6}$, then
\begin{equation}\label{expansion}
  \BARH_M(\eta)=\BARH_M^{(0)}(\eta)
  -\sum_{n\geq 1}\frac{1}{n!}\sum_{
  \genfrac{}{}{0cm}{2}{(R_1,\ldots,R_n)\in\Rr^n}{R_i\subset\LATTC}}
  \phi(R_1,\ldots,R_n)\prod_{i=1}^n \zeta(R_i)\COMMA
\end{equation}
where
\[
  \phi(R_1,\ldots,R_n)=\left\{\begin{array}{ll}1&{,n=1}\\
  {\sum_{g\in G_n}\prod_{\{i,j\}\in g}(1(R_i,R_j)-1)}&{,n>1}\end{array}\right.\PERIOD
\]
and $G_n$ is  the set of the generalized, connected graphs on $\{1,\ldots,n\}$ and
\[
  1(R_i,R_j)=\begin{cases} 0, & \{R_i\cap R_j\neq\emptyset\}\,;\\
                           1, & \{R_i\cap R_j=\emptyset\}\PERIOD
             \end{cases}
\]
\end{lemma}
\PROOF
The proof is the same as in \cite[Theorem 2]{Cammarota}
given \VIZ{perturbation} and the discussion preceding Lemma~\ref{lemma1}.
For later use we repeat the main steps of the proof. 
For the series \VIZ{expansion} we have
\begin{multline}
  \sum_{n\geq 1}\frac{1}{n!}\sum_{
  \genfrac{}{}{0cm}{2}{(R_1,\ldots,R_n)\in\Rr^n}{R_i\subset\LATTC}}
  \left|\phi(R_1,\ldots,R_n)\prod_{i=1}^n \zeta(R_i)\right| =\\
  =\sum_{R\in\Rr,R\subset\LATTC}|\zeta(R)|
  +\sum_{n\geq 2}\frac{1}{n!}\sum_{
  \genfrac{}{}{0cm}{2}{(R_1,\ldots,R_n)\in\Rr^n}{R_i\subset\LATTC}}
  \left|\phi(R_1,\ldots,R_n)\prod_{i=1}^n \zeta(R_i)\right|=:S\COMMA
\end{multline}
where for the first term ($n=1$) we have
\begin{equation}\label{n_1}
  \sum_{R\in\Rr,R\subset\LATTC}|\zeta(R)|
  \leq
  M\sum_{r\geq 1}\sup_{k\in\LATTC}\sum_{
  \genfrac{}{}{0cm}{2}{R\in\Rr,R\supset \{k\}}{|R|=r}}|\zeta(R)|
  \leq
  M(\delta+\sum_{r\geq 2}\delta^{r-1})\PERIOD
\end{equation}
For the second term ($n\geq 2$) using the estimate from 
\cite[p. 526, (3.6)]{Cammarota} we obtain
\begin{multline}
\sum_{n\geq 2}\frac{1}{n!}
\sum_{
\genfrac{}{}{0cm}{2}{(R_1,\ldots,R_n)\in\Rr^n}{R_i\subset\LATTC}
}
\left|\phi(R_1,\ldots,R_n)\prod_{i=1}^n \zeta(R_i)\right| \leq \\
\leq
\sum_{n\geq 2}\frac{1}{n!}
\sum_{
\genfrac{}{}{0cm}{2}{R\in\Rr}{R\subset\LATTC}
}
\sum_{
\genfrac{}{}{0cm}{2}{(R_1,\ldots,R_n)\in\Rr^n}{\exists R_i=R}
}
\left|\phi(R_1,\ldots,R_n)\prod_{i=1}^n \zeta(R_i)\right| \leq \\
\leq
\sum_{n\geq 2}\frac{n}{n!}
\sum_{
\genfrac{}{}{0cm}{2}{R\in\Rr}{R\subset\LATTC}
}
|\zeta(R)|
\sum_{(R_2,\ldots,R_n)\in\Rr^{n-1}}
\left|\phi(R_1,\ldots,R_n)\prod_{i=2}^n \zeta(R_i)\right| \leq \\
\leq
\sum_{n\geq 2}\frac{n}{n!}
\sum_{
\genfrac{}{}{0cm}{2}{R\in\Rr}{R\subset\LATTC}
}
|\zeta(R)|\frac{1}{2}\left(
2e\frac{5}{4}\frac{\delta}{1-\delta}\right)^{n-1}(n-2)!
|R|\EXP{|R|} \PERIOD
\end{multline}
Thus,
\begin{equation}\label{S}
S\leq M(\delta+\sum_{r\geq 2}\delta^{r-1})+M\frac{1}{2}
\sum_{n\geq 2}\frac{1}{n-1}
\left(
2e\frac{5}{4}\frac{\delta}{1-\delta}\right)^{n-1}
(\delta e+\sum_{r\geq 2}r(\delta e)^{r-1})\COMMA
\end{equation}
which concludes the proof.
\qed

\subsection{The small parameter of the cluster expansion}\label{errors}

In order to identify a small parameter we set 
\[
  E_{kl}(x-y):=J(x-y)-\bar{J}(k,l)\PERIOD
\]
The expansion \VIZ{star} suggests that the error
depends on both $\beta$ and
\[
  \sup_{k,l\in\LATTC}\sup_{x\in \CUBE_k,\, y\in\CUBE_l}|E_{kl}(x-y)|\PERIOD
\]

This quantity can be evaluated more explicitly for the special choice of
interaction potential defined in \VIZ{defJV1}.
\begin{lemma}\label{a-priori_error}
Assume that $J$ satisfies \VIZ{defJV1}--\VIZ{defJV2} then the coarse-grained
interaction potential $\bar J$ at the coarse-graining level $q$ approximates the potential $J$
with the error
\begin{equation}\label{errJb}
  |J(x-y) - \bar J(k,l)| \leq 2\frac{q}{L^{d+1}}
          \sup_
{\genfrac{}{}{0cm}{2}{x'\in \CUBE_k,}{y'\in \CUBE_l}}\|\nabla V(x'-y')\|
\end{equation}
for both $l\neq k$ and $l=k$.
\end{lemma}
\PROOF
Using the properties of the potential $V$, we expand
$V$ into the Taylor series, 
$$
V(z) = V(z') + (z-z') . \nabla V(z') + \BIGO(\|z-z'\|^2)\PERIOD
$$
Using the definition of $J$,  \VIZ{defJV1} and setting $z=x-y$ and $z'=x'-y'$, where
$x,x'\in\CUBE_k$ and $y,y'\in\CUBE_l$, we have

\begin{eqnarray*}
J(x-y) &=& \frac{1}{q^{2d}} \sum_{x'\in\CUBE_k}\sum_{y'\in\CUBE_l} J(x'-y') +\\
       && + \frac{1}{q^{2d}}\sum_{x'\in\CUBE_k}\sum_{y'\in\CUBE_l}
\frac{n}{L^d L} ((x-y) - (x'-y')).\nabla V(\frac{n}{L}(x'-y')) \\
       && + \frac{1}{q^{2d}}\sum_{x'\in\CUBE_k}\sum_{y'\in\CUBE_l}\BIGO\left(
\|\frac{n}{L}((x-y)-(x'-y'))\|^2\right)\COMMA
\end{eqnarray*}
and using the estimate
$\|(x-y) - (x'-y')\|\leq \|x-x'\| + \|y - y'\| \leq 2\max\{\mathrm{diam}\,(\CUBE_k)\}\sim\frac{2}{m}$
we obtain \VIZ{errJb} in both cases $l\neq k$ and $l=k$.
\qed

Our goal is to construct higher-order corrections to the Hamiltonian $\BARH_M^{(0)}$
based on the expansion \VIZ{expansion}.
First we check the condition \VIZ{condition} and identify
the parameter $\delta$. From the point of view of the cluster expansion theory 
we are in the  
high temperature, small activity regime 
because of the a priori estimate given by Lemma~\ref{a-priori_error}
on the coarse-grained approximation. Hence we expect
the small parameter $\delta$ to depend on the characteristics of the coarse-graining,
in particular on  the parameters $q$ and $L$.

\begin{lemma}[{\rm Identification of the small parameter $\delta$}]
\label{lemma2}
For $\zeta$ defined in \VIZ{activity}, condition
\VIZ{condition} holds with 
\[
\delta\sim\beta\sup_{k\in\LATTC}\sum_{l:\,l\neq k}|\sum_{x\in\CUBE_k,\,y\in\CUBE_l}E_{kl}(x-y)|
\]
In the particular case of Lemma~\ref{a-priori_error} we obtain
$\delta\sim\beta\frac{q^{d+1}}{L}\|\GRADV\|_\infty$,
where $q$ is the level of coarse-graining and
$L$ is the range of the interaction of the particles on the microscopic lattice.
\end{lemma}
\PROOF
To check the condition \VIZ{condition} we follow the analysis proposed
in \cite{PdLS}. 
We start with the $r=1$ case
\begin{multline*}
\sup_{k\in\LATTC}\sum_{
\genfrac{}{}{0cm}{2}{R\in\Rr,R\supset \{k\}}{|R|=1}}|\zeta(R)|
=\sup_{k\in\LATTC}\int |\EXP{-\beta\Delta_{kk}J(\sigma)}-1|\tilde{\rho}_k(d\sigma)\\
\sim\BIGO\LP\LP\frac{q}{L}\RP^d\beta\frac{q^{d+1}}{L}\|\GRADV\|_\infty\RP\COMMA
\end{multline*}
where we have used the fact that $|\EXP a-1|\leq |a|\EXP{|a|}$
and that $\EXP{\beta|\Delta_{kk}J(\sigma)|}\sim\BIGO(1)$.
For the general case $r\geq 2$ we have
\begin{multline}
\sup_{k\in\LATTC}\sum_{
\genfrac{}{}{0cm}{2}{R\in\Rr,R\supset \{k\}}{|R|=r}}
\left|\int\sum_{g\in G_R}\prod_{\{l_1,l_2\}\in g}f_{l_1l_2}(\sigma)
\prod_{\{k\}\in R}\tilde{\rho}_k(d\sigma)
\right|= \\
=
\sup_{k\in\LATTC}\frac{1}{(r-1)!}\sum_{
\genfrac{}{}{0cm}{2}{k_1=k,k_2,\ldots,k_r}{k_i\neq k_j,i\neq j}}
\left|\int
\sum_{g\in G_r}\prod_{\{i,j\}\in g}(\EXP{-\beta\Delta_{k_i k_j}J(\sigma)}-1)
\prod_{i=1}^r\tilde{\rho}_{k_i}(d\sigma)
\right|\COMMA \label{one}
\end{multline}
where we have assumed $R=\{k_1,\ldots,k_r\}$. We divide
by $(r-1)!$ because we can relabel $k_2,\ldots,k_r$.  Note that by
$G_r$ we denote the set of all connected graphs in $\{1,2\ldots,r\}$
which may have self-loops $\{j,j\}$ for any combination of
$j=1,2\ldots,r$ (from now on we simply call them loops instead of
self-loops).  In the sequel we need the following estimate 
due to Lemma~\ref{a-priori_error} we have that for a
fixed $k\in\LATTC$
\[
\sum_{l:\,l\neq k}|\Delta_{kl}J(\sigma)|
  \sim \LP\frac{L}{q}\RP^d q^{2d}\frac{1}{L^d}\BIGO\LP\frac{q}{L}\RP
  \sim \BIGO\LP\frac{q^{d+1}}{L}\|\GRADV\|_\infty\RP\COMMA
\]
which implies that
\begin{equation}\label{C}
\sup_{k\in\LATTC}\sum_{l:\,l\neq k}|\Delta_{kl}J(\sigma)|\leq C\COMMA
\end{equation}
for every $\sigma$ and with 
$C\sim\sup_{k\in\LATTC}\sum_{l:\,l\neq k}|\sum_{x\in\CUBE_k,\,y\in\CUBE_l}E_{kl}(x-y)|$ 
or in the special case of $J$ given in terms of $V$ 
$C\sim\frac{q^{d+1}}{L}\|\GRADV\|_\infty$.

Now we denote by $G^0_r$ the set of graphs $g^0$
without loops and we describe by 
$g\sim g^0$ the corresponding graph $g$ with loops such that
$g=g^0$ if we remove all loops $\{j,j\}$ from $g$.
We have
\begin{multline}
\sum_{g\in G_r}\prod_{\{i,j\}\in g}(\EXP{-\beta\Delta_{k_i k_j}J(\sigma)}-1)= \\
=\sum_{g^0\in G^0_r}\prod_{\{i,j\}\in g^0}(\EXP{-\beta\Delta_{k_i k_j}J(\sigma)}-1)
\left(1+\sum_{g\sim g^0}\prod_{\{j,j\}\in g/ g^0}
(\EXP{-\beta\Delta_{k_j k_j}J(\sigma)}-1)
\right)\PERIOD
\end{multline}
We observe that for any $g^0\in G_r^0$ with $\supp(g^0)=r$ the terms in the 
parenthesis are  bounded by $2^r$, since
for every $g\sim g^0$ we have that 
$\prod_{\{i,j\}\in g/ g^0}\beta|\Delta_{k_i k_j}J(\sigma)|\leq 1$
(after the choice we will make at the end of the proof)
and the summation $\sum_{g\sim g^0}$ has 
$\sum_{k=1}^r \binom{r}{k}=2^r-1$ many terms.
Then, for the sum $
\sum_{g^0\in G^0_r}\prod_{\{i,j\}\in g^0}(\EXP{-\beta\Delta_{k_i k_j}J(\sigma)}-1)
$
we can use the tree-graph equality (e.g., \cite{PdLS}[p.5, Theorem 2]) to estimate the sum
over connected graphs  with respect to
the sum over the corresponding maximal spanning trees.
From \cite[Theorem 2]{PdLS} we have that
\begin{multline}
  \sum_{g^0\in G^0_r}\prod_{\{i,j\}\in g^0}(\EXP{-\beta\Delta_{k_i k_j}J(\sigma)}-1)= \\
 =\sum_{\tau^0\in\Tt^0_r}\prod_{\{i,j\}\in\tau^0}(-\beta\Delta_{k_i k_j}J(\sigma))\int_0^1 \!\!dt_1
  \ldots\int_0^1 \!\!dt_{r-1}
\!\!\!\!\!\!\sum_{
\genfrac{}{}{0cm}{2}{X_1,\ldots,X_{r-1}}{\mbox{\scriptsize comp. with} \,\tau^0}
} \\
t_1^{b_1-1}\ldots t_{r-1}^{b_{r-1}-1}\EXP{-W(X_1,\ldots,X_{r-1},t_1,\ldots,t_{r-1})}\COMMA\label{tree-graph}
\end{multline}
where $g^0$ is a connected graph in $\{1,2\ldots,r\}$ without loops and 
$\tau^0$ is the corresponding tree graph.
The summation is over all sequences $X_1,\ldots,X_{r-1}$ of increasing subsets of $\{1,2\ldots,r\}$
such that $X_1={1},\,X_n\subset X_{n+1}$ and $|X_n|=n$ that are compatible with $\tau^0$.
A sequence $X_1,\ldots,X_{r-1}$ is compatible with a given tree $\tau^0$ if, 
for all $n=1,\ldots,r-1$, $X_n$ contains exactly $n-1$ links of $\tau^0$.
We also say that a link $\{i,j\}$ crosses $X_n$ if $i\in X_n$ and $j\notin X_n$
or vice versa; then $b_n$ is the number of links of $\tau^0$ which cross $X_n$, and
\[
W(X_1,\ldots,X_{r-1},t_1,\ldots,t_{r-1})=\sum_{1\leq i< j\leq r}
t_1(\{i,j\})\ldots t_{r-1}(\{i,j\})\beta\Delta_{k_i k_j}J(\sigma)\COMMA
\]
where
\[
t_n(\{i,j\})=\left\{\begin{array}{ll} {t_n\in[0,1]} 
&{,\mbox{if}\, \{i,j\}\,\mbox{crosses}\, X_n}\\
1&{,\mbox{otherwise}}\end{array}\right.\PERIOD
\]
Then from \VIZ{C} we obtain that
\begin{eqnarray*}
W(X_1,\ldots,X_{r-1},t_1,\ldots,t_{r-1}) & \leq &
\sum_{1\leq k< l\leq r}|\beta\Delta_{kl}J(\sigma)|\\
& \leq & r \beta\sup_{k\in\LATTC}\sum_{l=1, l\neq k}^r |\Delta_{kl}J(\sigma)|
\sim r\beta C\COMMA
\end{eqnarray*}
where $C\sim\BIGO(\frac{q^{d+1}}{L}\|\GRADV\|_{\infty})$. 
Thus, uniformly in $X_1,\ldots,X_{r-1},t_1,\ldots,t_{r-1}$ we have that
\[
\EXP{-W(X_1,\ldots,X_{r-1},t_1,\ldots,t_{r-1})}\leq \EXP{r\beta C}\PERIOD
\]
Furthermore, from \cite[Lemma 4]{PdLS} we have that
\[
\int_0^1 dt_1
\ldots\int_0^1 dt_{r-1}\!\!\!\!\!\!\!\!
\sum_{
\genfrac{}{}{0cm}{2}{X_1,\ldots,X_{r-1}}{{\scriptstyle\mbox{\scriptsize comp. with}} \,\tau^0}
}
t_1^{b_1-1}\ldots t_{r-1}^{b_{r-1}-1}= 1\PERIOD
\]
This leads to the following tree-graph inequality for the case of
graphs with loops
\begin{eqnarray*}
\left|\sum_{g\in G_r}\prod_{\{i,j\}\in g}(\EXP{-\beta\Delta_{k_i k_j}J(\sigma)}-1)\right|
& \leq & 2^r
\left|\sum_{g^0\in G^0_r}\prod_{\{i,j\}\in g^0}(\EXP{-\beta\Delta_{k_i k_j}J(\sigma)}-1)\right|\\
& \leq  & \mbox{} 2^r\EXP{r\beta C}\beta\sum_{\tau^0\in\Tt^0_r}\prod_{\{i,j\}\in
\tau^0}|\Delta_{k_i k_j}J(\sigma)|\PERIOD
\end{eqnarray*}
From \VIZ{C} we have that for any tree $\tau^0$
\[
\sup_{k\in\LATTC}\sum_{
\genfrac{}{}{0cm}{2}{k_1=k,k_2,\ldots,k_r}{k_i\neq k_j,i\neq j}
}
\prod_{\{i,j\}\in\tau^0}|\Delta_{k_i,k_j}J(\sigma)|\leq C^{r-1}\PERIOD
\]
Using the Cayley formula $\sum_{\tau^0\in\Tt_r^0}1=r^{r-2}$,
we conclude that
\[
\sup_{k\in\LATTC}\sum_{
\genfrac{}{}{0cm}{2}{R\in\Rr,R\supset \{k\}}{|R|=r}
}|\zeta(R)|\leq 
2^r\frac{\EXP{r\beta C}}{(r-1)!}r^{r-2}(\beta C)^{r-1}
\PERIOD
\]
Thus, we require
\[
2\beta C\EXP{\beta C+1}
\leq\delta\COMMA
\]
so it suffices to consider
$\delta\sim\beta \frac{q^{d+1}}{L}\|\GRADV\|_\infty$.
\qed 

%
%
\subsection{Effective Hamiltonians and specific entropy estimates}
\label{strategy}
In this section we use the results of Lemma~\ref{lemma1} and \ref{lemma2}
and expand the Hamiltonian $\bar{H}_M(\eta)$ around $\bar{H}_M^{(0)}$ in
powers of $\delta$ by truncating the series \VIZ{expansion} at the levels
of the number of interacting clusters and of the length of the clusters.
We obtain the series
$$
\bar{H}_M(\eta) =  \bar{H}_M^{(0)} + \bar{H}_M^{(0)}  + \cdots +
\bar{H}_M^{(p)} + MO(\delta)\COMMA
$$
and define the corresponding Gibbs measures by
\begin{equation}\label{cg_one}
  \BARIT{\mu}_{M,\beta}^{(p)}(d\eta)=\frac{1}{\BARIT{Z}^{(p)}_M}
  \EXP{-\beta(\BARH_M^{(0)}(\eta)+
  \ldots+\BARH_M^{(p)}(\eta))}\
  \BARIT{P}_M(d\eta)\COMMA
\end{equation}
where $p=1,2, \cdots$ and  $\BARIT{Z}_M^{(p)}$ is the corresponding partition
function. In the following proposition we derive explicit formulas up to the
order $\BIGO(\delta^3)$ and state the corresponding error estimates for the
relative entropy.  Higher order corrections can be computed in a
straightforward although tedious way. 
For the sake of simplicity we restrict ourselves to the case $d=1$, 
however similar computations and results are derived in any dimension.
\begin{proposition}
[{\rm Corrections to $\BARH_M^{(0)}$}]
\label{proposition}
For $\BARH_M^{(0)}$ we have the error estimate
\begin{equation}\label{H_0}
\BARH_M(\eta)=\BARH_M^{(0)}(\eta)+M\BIGO(\delta)\PERIOD
\end{equation}
Furthermore, the first and second order corrections to the Hamiltonian 
$\BARH_M^{(0)}$ are given by
\begin{eqnarray}\label{second}
\BARH_M(\eta) & = & \BARH_M^{(0)}(\eta)+\BARH_M^{(1)}(\eta)+M\BIGO(\delta^2)\nonumber\\
\BARH_M(\eta) & = & \BARH_M^{(0)}(\eta)+\BARH_M^{(1)}(\eta)+\BARH_M^{(2)}(\eta)
+M\BIGO(\delta^3)
\COMMA
\end{eqnarray}
where
\begin{eqnarray}\label{H_1}
\BARH_M^{(1)}(\eta) & = &
        -\sum_{k\in\LATTC}\int f_{kk}(\sigma)\tilde{\rho}_k(d\sigma)-{}\nonumber\\
& &   {}-\sum_{k}\!\!\!\sum_{l\in B^+_{\frac{L}{q}}(k)}\int \! 
              (f_{kl}+\! f_{kl}f_{kk}+\! f_{kl}f_{ll}+\! f_{kl}f_{kk}f_{ll})
              \tilde{\rho}_k(d\sigma)\tilde{\rho}_l(d\sigma)
\end{eqnarray}
and
\begin{eqnarray}
\BARH_M^{(2)}(\eta)&=& 
       {}\frac{1}{2}\sum_{k\in\LATTC}\left(\int f_{kk}(\sigma)\tilde{\rho}_k(d\sigma)
       \right)^2 +{} \nonumber\\
   & &{}+\frac{1}{2}\sum_{k}\sum_{l\in B^+_{\frac{L}{q}}(k)}
        \int f_{kk}(\sigma)\tilde{\rho}_k(d\sigma)\int f_{kl}(\sigma)
        \tilde{\rho}_k(d\sigma)\tilde{\rho}_l(d\sigma)+{}\nonumber\\
   & &{}+\frac{1}{2}\sum_{k}\sum_{l_1\in B^+_{\frac{L}{q}}(k)}
       \sum_{l_2\in B^+_{\frac{L}{q}}(k)} \nonumber \\
   & & \;\;\;\;\int (f_{kl_1}+f_{kk}f_{kl_1}+f_{kl_1}f_{l_1l_1}+f_{kk}f_{kl_1}f_{l_1l_1})
           \tilde{\rho}_k(d\sigma)\tilde{\rho}_{l_1}(d\sigma)\times\nonumber \\
   & & \;\;\;\;\times
       \int (f_{kl_2}+f_{kk}f_{kl_2}+f_{kl_2}f_{l_2l_2}+f_{kk}f_{kl_2}f_{l_2l_2})
       \tilde{\rho}_k(d\sigma)\tilde{\rho}_{l_2}(d\sigma) +{}\nonumber \\
   & &{}+\sum_{k}\sum_{l_1\in B^+_{\frac{L}{q}}(k)}
       \sum_{l_2\in B_{\frac{L}{q}}(l_1)\cup B^+_{\frac{L}{q}}(k)} \nonumber \\
   & & \;\;\int(f_{kl_1}f_{l_1l_2}+f_{kl_1}f_{kl_2}+f_{kl_2}f_{l_1l_2}+[\ldots]) 
       \tilde{\rho}_k(d\sigma)\tilde{\rho}_{l_1}(d\sigma)\tilde{\rho}_{l_2}(d\sigma)\nonumber \\
   &=& I_1+I_2+I_3+I_4 \COMMA \label{H_2}
\end{eqnarray}
where $[\ldots]$ means the previous three terms with all possible combinations of loops.
The corresponding relative entropy error is
\[
\RELENTR(\BARIT{\mu}_{M,\beta}^{(p)}|\mu_{N,\beta}\circ \COP^{-1})
\sim\BIGO\LP\frac{1}{q}\delta^{p+1}\RP\PERIOD
\]
\end{proposition}
\PROOF
The error in \VIZ{H_0} is given by \VIZ{S}.
To get an error $M\BIGO(\delta^2)$ we need to include 
the terms $n=1$, $r=1,2$, i.e.,
\[
\BARH_M^{(1)}(\eta)=-\sum_{R:\,|R|=1,2}\zeta(R)\PERIOD
\]
Similarly, for the $M\BIGO(\delta^3)$ error we need to include \VIZ{H_1} 
together with the terms $n=2$, $r=1,2$ and $n=1$, $r=3$
\[
\BARH_M^{(2)}(\eta)=-\frac{1}{2}\sum_{
\genfrac{}{}{0cm}{2}{R_1,R_2}{|R_i|=1,2,i=1,2}}
\phi(R_1,R_2)\zeta(R_1)\zeta(R_2)
-\sum_{R:\,|R|=3}\zeta(R)\COMMA
\]
where $\phi(R_1,R_2)=0$ if $R_1\cap R_2=\emptyset$
and $-1$ if $R_1\cap R_2\neq\emptyset$,
which occurs when $R_1=R_2=\{k\}$ for $k\in\LATTC$, in the case of $|R_i|=1$,
$i=1,2$. If $|R_i|=2$, $i=1,2$ we have that $R_1=\{k,l_1\}$ and $R_1=\{k,l_2\}$.
For the relative entropy error
for $\BARH_M^{(p)}$ with $p=0,1,\ldots$, we have
\begin{eqnarray}\label{relative_entropy}
\RELENTR(\BARIT{\mu}_{M,\beta}^{(p)}|\mu_{N,\beta}\circ \COP^{-1}) & = & 
\frac{1}{N}\int\log\frac{\frac{1}{\BARIT{Z}_M^{(p)}}\EXP{-\beta\BARH_M}\rho(\eta)}
{\frac{1}{Z_N}\int_{\{\COP(\sigma)=\eta\}}\EXP{-\beta H_N(\sigma)}P_N(d\sigma)}
\BARIT{\mu}_{M,\beta}^{(p)}(d\eta)\nonumber\\
& = & \frac{1}{N}\int\log\frac{\frac{1}{\BARIT{Z}_M^{(p)}}\EXP{-\beta\BARH_M^{(p)}}\rho(\eta)}
{\frac{1}{\BARIT{Z}_M}\EXP{-\beta\BARH_M}\rho(\eta)}
\BARIT{\mu}_{M,\beta}^{(p)}(d\eta)\nonumber\\
& = & 
\frac{1}{N}\log\frac{\BARIT{Z}_M}{\BARIT{Z}_M^{(p)}}
+\frac{1}{N}\EXPECT_{\BARIT{\mu}^{(p)}_{M,\beta}}[\beta(\BARH_M-\BARH_M^{(p)})]\COMMA
\end{eqnarray}
where for the partition functions we have
\begin{eqnarray}\label{partitions}
\BARIT{Z}_M&=&
\sum_{\eta} \EXP{-\beta\BARH_M}\rho(\eta)=
\BARIT{Z}_M^{(p)}\sum_{\eta}
\frac{1}{\BARIT{Z}_M^{(p)}}\EXP{-\beta\BARH_M^{(p)}}\EXP{-\beta(\BARH-\BARH_M^{(p)})}\rho(\eta)=
    \nonumber \\
&=& \BARIT{Z}_M^{(p)}\,\EXPECT_{\BARIT{\mu}^{(p)}_{M,\beta}} [\EXP{-\beta(\BARH_M-\BARH_M^{(p)})}]
\PERIOD
\end{eqnarray}
Thus, for the different choices of $p$ in $\BARH_M^{(p)}$,
given the fact that the estimates \VIZ{H_0} and \VIZ{second} are uniform in $\eta$ we get that
\begin{eqnarray*}
&& \log\frac{\BARIT{Z}_M}{\BARIT{Z}_M^{(p)}}\sim\log\EXPECT_{\BARIT{\mu}^{(p)}_{M,\beta}}\left[
       \EXP{M\BIGO(\delta^{p+1})}\right]\sim M\BIGO(\delta^{p+1})\COMMA\\
&& \EXPECT_{\BARIT{\mu}^{(p)}_{M,\beta}}[\BARH_M-\BARH_M^{(p)}]\sim M\BIGO(\delta^{p+1})\COMMA
\end{eqnarray*}
where $p=0,1,2$.
\qed

%
%
\section{Numerical schemes for coarse-graining and a posteriori estimates}\label{schemes}
We conclude the proof of Theorem~\ref{mainresult}
by giving explicit formulas for the correction terms  \VIZ{H_1} and \VIZ{H_2}
and estimating them with respect to the small parameter $\epsilon$.
The corrections
consist of combinations of $f_{kl}=\EXP{-\beta\Delta_{kl}J(\sigma)}-1$
for all  $k$ and $l$.
Since the exponent $-\beta\Delta_{kl}J(\sigma)$ 
in \VIZ{f} is small we have that
\begin{equation}\label{p}
   \EXP{-\beta\Delta_{kl}J(\sigma)}-1=\sum_{p=1}^{\infty}
   \frac{1}{p!}(-\beta\Delta_{kl}J(\sigma))^p
   \,\,\mbox{with}\,\,\,
   \Delta_{kl}J(\sigma)\sim
   \BIGO(q^{2d}\frac{q}{L^{d+1}}\|\GRADV\|_\infty)\PERIOD
\end{equation}
The key point of our calculations relies on the cancellation
\begin{equation}\label{deltaJiszero}
    \int\Delta_{kl}J(\sigma)\tilde{\rho}_k(d\sigma)\tilde{\rho}_l(d\sigma)=0\COMMA
\end{equation}
following from the definition of $\BARH^{(0)}_M$ in \VIZ{cg_Hamiltonian}.

In order to keep expressions simple we express our results in terms of the
variables
$\alpha(k)$, the number of spins $\sigma(x)=1$ in the coarse cell $\CUBE_k$), and 
$\vita(k)= q-\alpha(k)$.
Next we define the following quantities
\begin{eqnarray}
E_1(\alfa)&:=&\EXPECT[\sigma(x)|\alfa]=\frac{2\alfa-q}{q} \label{e1} \\
E_2(\alfa)&:=&\EXPECT[\sigma(x)\sigma(y)|\alfa]=
              \frac{\alfa(\alfa-1)-2\alfa\vita+\vita(\vita-1)}{q(q-1)} \label{e2} \\
E_3(\alfa)&:=&\EXPECT[\sigma(x)\sigma(y)\sigma(z)|\alfa]= \nonumber \\
          & =&\frac{\alfa(\alfa-1)(\alfa-2)-3\alfa(\alfa-1)\vita+3\alfa(\vita-1)\vita-(\vita-2)(\vita-1)\vita}
{q(q-1)(q-2)} \label{e3}\\
E_4(\alfa)&:=&\EXPECT[\sigma(x)\sigma(y)\sigma(z)|\alfa]\nonumber = \\
          & =& \frac{\alfa(\alfa-1)(\alfa-2)(\alfa-3)-4\alfa(\alfa-1)(\alfa-2)\vita}
               {q(q-1)(q-2)(q-3)}+{} \nonumber\\
          & +&\frac{6\alfa(\alfa-1)(\vita-1)\vita-4\alfa(\vita-2)(\vita-1)\vita
               +\vita(\vita-1)(\vita-2)(\vita-3)}{q(q-1)(q-2)(q-3)}\label{e4}
\end{eqnarray}
which are all of order $\BIGO(1)$.
Furthermore, we introduce the notation
\begin{eqnarray}
j^1_{kl}&:=&\sum_{\genfrac{}{}{0cm}{2}{x\in\CUBE_k}{y\in\CUBE_l}}
             (J(x-y)-\bar{J}(k,l))^2 \label{j1} \\
j^2_{kl}&:=&\sum_{\genfrac{}{}{0cm}{2}{x\in\CUBE_k}{y,y'\in\CUBE_l}}
             (J(x-y)-\bar{J}(k,l))(J(x-y')-\bar{J}(k,l)) \label{j2} \\
j^2_{k_1k_2k_3}&:=&\sum_{\genfrac{}{}{0cm}{2}{x\in\CUBE_{k_1}}{y\in\CUBE_{k_2},z\in\CUBE_{k_3}}}
                 (J(x-y)-\bar{J}(k_1,k_2))(J(y-z)-\bar{J}(k_2,k_3)) \label{j23}
\end{eqnarray}
If $k_1=k_2$ then we also impose that for $x,y\in\CUBE_{k_1}$ we
have $y\neq x$.

Note that these quantities have various symmetries, for example, $j^2_{lk}=j^2_{kl}$ 
or $j^1_{kl} = \tilde j^1_{k-l} = \tilde j^1_{l-k}$ for some function $\tilde j^1$ and
similarly $j^2_{kl}$ depends also only on $|k-l|$, 
moreover $j^2_{k_1k_2k_3} = \tilde j^2_{k_1-k_2,k_3-k_2}$.

Based on \VIZ{errJb1} we get first estimates in terms of $q$ and $L$ 
\begin{align*}
j^1_{kl}       & \sim \LP\frac{q}{L}\RP^{2d+2}, &
j^2_{kl}       & \sim q^{3d}\LP\frac{q}{L^{d+1}}\RP^2\sim q^d\LP\frac{q}{L}\RP^{2d+2}, &
j^2_{k_1k_2k_3}& \sim q^d\LP\frac{q}{L}\RP^{2d+2}\PERIOD
\end{align*}

\noindent{\bf Proof of Theorem~\ref{mainresult}:}
Starting from the formulas for  $\BARH_M^{(1)}$ and $\BARH_M^{(2)}$
in Proposition~\ref{proposition} 
we expand the factors $f_{kl}$ and re-estimate the corrections using the
cancellations in \VIZ{deltaJiszero}.
Recalling that $\epsilon\sim \frac{q}{L}\beta\|\GRADV\|_{\infty}$,
we have
\begin{eqnarray*}
&\mbox{(i) }&\sum_{k\in\LATTC}\int f_{kk}(\sigma)\tilde{\rho}_k(d\sigma)
  \sim M \LP q^{2d}\frac{q}{L^{d+1}}\beta\|\GRADV\|_{\infty}\RP^2
  \sim M q^{2d}\LP\frac{q}{L}\RP^{2d}\epsilon^2
\\
&\mbox{(ii) }&\sum_{k}\sum_{l\in B_{\frac{L}{q}}(k)}\int (f_{kl}+f_{kl}f_{kk}
   +f_{kl}f_{ll})\tilde{\rho}_k(d\sigma)\tilde{\rho}_l(d\sigma) \sim \\
&& \sim M \LP\frac{L}{q}\RP^d \LP q^{2d}\frac{q}{L^{d+1}}\beta\|\GRADV\|_{\infty}\RP^2 
   \sim M q^{2d}\LP\frac{q}{L}\RP^d \epsilon^2
\end{eqnarray*}
while the term involving $f_{kl}f_{kk}f_{ll}$ is of a higher order.

Therefore it suffices to include only the first two contributions
which we also compute explicitly.
With a slight abuse of notation we still denote the new terms by $\BARH_M^{(1)}$.
Collecting all the terms gives us
\begin{multline}
-\BARH_M^{(1)}=\beta\sum_k\int\frac{1}{8}(\Delta_{kk}J(\sigma))^2\tilde{\rho}_k(d\sigma)
+\beta\sum_{k<l}\int\frac{1}{2}(\Delta_{kl}J(\sigma))^2\tilde{\rho}_k(d\sigma)\tilde{\rho}_l(d\sigma) +{}\\
   {}+\beta\sum_{k<l}\int\frac{1}{2}\Delta_{kk}J(\sigma)\Delta_{kl}J(\sigma)\tilde{\rho}_k(d\sigma)
      \tilde{\rho}_l(d\sigma)+{}\\
   {}+\beta\sum_{k<l}\int\Delta_{kl}J(\sigma)\frac{1}{2}\Delta_{ll}J(\sigma)\tilde{\rho}_k(d\sigma)
      \tilde{\rho}_l(d\sigma)\PERIOD
\end{multline}
Explicit calculations of these integrals are detailed in Appendix~A, yielding
\begin{eqnarray*}
-\BARH_M^{(1)}
     &=&\frac{\beta}{8}\sum_k4j_{kk}^2\left[-E_4(\alfa(k))+E_2(\alfa(k))\right]
                  +2j_{kk}^1\left[E_4(\alfa(k))+1-2E_2(\alfa(k))\right]+{} \nonumber\\
     & &
\begin{aligned}
    {}+\frac{\beta}{2}\sum_{k<l}& j_{kl}^1\left[E_2(\alfa(k))E_2(\alfa(l))-E_2(\alfa(l))-E_2(\alfa(k))+1)\right]+{} \\
     &{}+j_{kl}^2\left[-2E_2(\alfa(k))E_2(\alfa(l))+E_2(\alfa(k))+E_2(\alfa(l))\right]+{}
\end{aligned}\nonumber\\
     & &{}+\frac{\beta}{2}\sum_{k<l}2j_{kkl}^2\left[-E_3(\alfa(k))E_1(\alfa(l))
                                                    +E_1(\alfa(k))E_1(\alfa(l))\right]+{}\nonumber\\
     & &{}+\frac{\beta}{2}\sum_{k<l}2j_{llk}^2\left[-E_3(\alfa(l))E_1(\alfa(k))+
                                                     E_1(\alfa(l))E_1(\alfa(k))\right]\PERIOD
\end{eqnarray*}
Since $j_{kkl}^2=j_{llk}^2$ this expression further simplifies to
\begin{eqnarray}
-\BARH_M^{(1)}
     &=&\frac{\beta}{8}\sum_k 4j_{kk}^2\left[-E_4(\alfa(k))+E_2(\alfa(k))\right]
                   +2j_{kk}^1\left[E_4(\alfa(k))+1-2E_2(\alfa(k))\right]+{}\nonumber \\
     & &
\begin{aligned}
   {}+\frac{\beta}{2}\sum_{k<l}&  j_{kl}^1\left[E_2(\alfa(k))E_2(\alfa(l))-E_2(\alfa(l))-E_2(\alfa(k))+1\right]+{}\nonumber \\
              &+{}j_{kl}^2\left[-2E_2(\alfa(k))E_2(\alfa(l))+E_2(\alfa(k))+E_2(\alfa(l))\right]+{}
\end{aligned} \\
     & &
\begin{aligned}
     {}+\frac{\beta}{2}\sum_{k,l\neq k} j_{kkl}^2 &\left[-E_3(\alfa(k))E_1(\alfa(l))
                                  +2E_1(\alfa(k))E_1(\alfa(l))-{}\right.\nonumber \\
          &{}-\left. E_3(\alfa(l))E_1(\alfa(k))\right]\PERIOD
\end{aligned}
\end{eqnarray}
By the estimates on the terms $j^2_{kl}$ and the fact that
the terms $E_i(\alfa)$ are of order one, and 
by counting the summations we obtain
\begin{eqnarray*}
\BARH_M^{(1)} & \sim & M q^{3d}\LP\frac{q}{L^{d+1}}\beta\|\GRADV\|_{\infty}\RP^2
+ M \LP\frac{L}{q}\RP^d q^{3d}\LP\frac{q}{L^{d+1}}\beta\|\GRADV\|_{\infty}\RP^2 \sim \\
& \sim & M q^d \LP\frac{q}{L}\RP^{2d}\epsilon^2
+ M q^d \LP\frac{L}{q}\RP^d \LP\frac{q}{L}\RP^{2d}\epsilon^2
\sim N\LP\frac{q}{L}\RP^d\epsilon^2\COMMA
\end{eqnarray*}
hence we have gained an extra $q^d$ as compared with the previous estimate.
In the same spirit we estimate the terms of the correction $\BARH_M^{(2)}$. From \VIZ{H_2}
we have
\begin{align*}
& I_1 \sim  M \left(q^{2d}\frac{q}{L^{d+1}}\beta\|\GRADV\|_{\infty}\right)^4\COMMA \;\;
& I_2 \sim  M \left(\frac{L}{q}\right)^d\left(q^{2d}\frac{q}{L^{d+1}}\beta\|\GRADV\|_{\infty}\right)^4\COMMA \\
& I_3 \sim  M \left(\frac{L}{q}\right)^{2d}\left(q^{2d}\frac{q}{L^{d+1}}\beta\|\GRADV\|_{\infty}\right)^4\COMMA \;\;
& I_4 \sim  M \left(\frac{L}{q}\right)^{2d}\left(q^{2d}\frac{q}{L^{d+1}}\beta\|\GRADV\|_{\infty}\right)^2\PERIOD 
\end{align*}
Obviously, for $\BARH_M^{(2)}$ it suffices to consider only the
term $I_4$.
Abusing slightly the  notation we call this term $\BARH_M^{(2)}$  and
we calculate it explicitly
\begin{multline}
\BARH_M^{(2)}  =  \beta\sum_{k_1}\sum_{k_2>k_1}\sum_{k_3>k_2}
       \int \Delta_{k_1k_2}J(\sigma)\Delta_{k_2k_3}J(\sigma)\tilde{\rho}_{k_1}(d\sigma)
         \tilde{\rho}_{k_2}(d\sigma)\tilde{\rho}_{k_3}(d\sigma)+{} \nonumber \\
       +\int \Delta_{k_2 k_3}J(\sigma)\Delta_{k_3 k_1}J(\sigma)\tilde{\rho}_{k_2}(d\sigma)
          \tilde{\rho}_{k_3}(d\sigma)\tilde{\rho}_{k_1}(d\sigma)+{} \nonumber \\
       +\int \Delta_{k_3 k_1}J(\sigma)\Delta_{k_1 k_2}J(\sigma)\tilde{\rho}_{k_3}(d\sigma)
\tilde{\rho}_{k_1}(d\sigma)\tilde{\rho}_{k_2}(d\sigma)
\end{multline}
Similarly, as for the previous term, we evaluate the integrals in terms of $E_k$ and obtain
$$
\begin{aligned}
\BARH_M^{(2)} =  \beta\sum_{k_1}&\sum_{k_2>k_1}\sum_{k_3>k_2}  \\ 
          & j^2_{k_1 k_2 k_3}(-E_1(\alfa(k_1))E_2(\alfa(k_2))E_1(\alfa(k_3))
              +E_1(\alfa(k_1))E_1(\alfa(k_3))) \nonumber \\
          +{}& j^2_{k_2 k_3 k_1}(-E_1(\alfa(k_2))E_2(\alfa(k_3))E_1(\alfa(k_1))
              +E_1(\alfa(k_2))E_1(\alfa(k_1))) \\
          +{}& j^2_{k_3 k_1 k_2}(-E_1(\alfa(k_3))E_2(\alfa(k_1))E_1(\alfa(k_2))
              +E_1(\alfa(k_3))E_1(\alfa(k_2)))\PERIOD
\end{aligned}
$$
Using the estimates above we obtain the order of the correction
\[
\BARH_M^{(2)}\sim 
M \LP\frac{L}{q}\RP^{2d} q^{3d}\LP\frac{q}{L^{d+1}}\beta\|\GRADV\|_{\infty}\RP^2
\sim M q^d \epsilon^2\PERIOD
\]

Thus, from the above we deduce two results:
the first one that due to the fact that the leading order terms 
that participate in the $\BIGO(\delta)$ error vanish, the choice
of the initial Hamiltonian
$\BARH_M^{(0)}$ yields a second-order accurate approximation, in other words
\[
\frac{\beta}{N}(\BARH_M-\BARH_M^{(0)})\sim\BIGO(\epsilon^2)\PERIOD
\]
The second one is that to obtain a $\BIGO(\epsilon^3)$ error, the
corresponding correction term have to include both
$\BARH_M^{(1)}$ and $\BARH_M^{(2)}$.
In such a case we obtain
\[
\frac{\beta}{N}\LP\BARH_M-(\BARH_M^{(0)}+\BARH_M^{(1)}+\BARH_M^{(2)})\RP\sim\BIGO(\epsilon^3)\PERIOD
\]
In the same spirit we can derive error estimates for the additional  terms
in the expansion.
\qed

In Section~\ref{sims} we discuss the computational complexity of higher-order terms.
Here we present only two schemes that are relevant for practical implementations. 
%
\begin{scheme}[{\rm  2nd-order coarse-graining}]\label{scheme1}
The 2nd-order coarse-graining algorithm has the following characteristics
\begin{enumerate}
\item Hamiltonian: $\BARH^{(0)}_M$, given by \VIZ{cg_Hamiltonian}.
\item Gibbs measure: $\BARIT{\mu}_{M,\beta}^{(0)}$, given by \VIZ{cg_one} for $p=1$.
\item Relative entropy error
\[
\RELENTR(\BARIT{\mu}_{M,\beta}^{(0)}|\mu_{N,\beta}\circ \COP^{-1})\sim 
\BIGO(\epsilon^2)\PERIOD
\]
\end{enumerate}
\end{scheme}
Scheme~\ref{scheme1} is the coarse-graining algorithm that has been extensively 
studied in \cite{KMVPNAS,KMV,kt,KPS}.
The novelty presented in this paper is the sharper error estimate that shows
that the error is of the order $\BIGO(\epsilon^2)$; this estimate readily follows 
from the cancellation due to \VIZ{deltaJiszero} as is seen from the calculations above.
%
\begin{scheme}[{\rm 3rd-order coarse-graining}]\label{scheme2}
We construct a higher-order Monte Carlo algorithm with the following characteristics
\begin{enumerate}
\item Hamiltonian: $\BARH_M^{(0)}+\BARH_M^{(1)}+\BARH_M^{(2)}$\COMMA
where the corrections are 
\begin{eqnarray}\label{cor1}
-\BARH_M^{(1)}(\eta)
     &=&
\begin{aligned}[t]
   \frac{\beta}{8}\sum_k &4 j_{kk}^2\left[-E_4(\alfa(k))+E_2(\alfa(k))+\right] \\
     +{}&                 2 j_{kk}^1\left[E_4(\alfa(k))+1-2E_2(\alfa(k))\right]+{}
\end{aligned} \\
     & &
\begin{aligned}
   {}+\frac{\beta}{2}\sum_{k<l}&j_{kl}^1\left[E_2(\alfa(k))E_2(\alfa(l))-E_2(\alfa(l))-E_2(\alfa(k))+1\right]+{}\nonumber \\
              +{}& j_{kl}^2\left[-2E_2(\alfa(k))E_2(\alfa(l))+E_2(\alfa(k))+E_2(\alfa(l))\right]+{}
\end{aligned} \\
     & &
\begin{aligned}
     {}+\frac{\beta}{2}\sum_{k,l\neq k} j_{kkl}^2 &\left[-E_3(\alfa(k))E_1(\alfa(l))
                                  +2E_1(\alfa(k))E_1(\alfa(l))-{}\right.\nonumber \\
          &{}-\left. E_3(\alfa(l))E_1(\alfa(k))\right]\PERIOD
\end{aligned}
\end{eqnarray}
and
\begin{equation}\label{cor2}
\begin{aligned}
& \BARH_M^{(2)}(\eta) =  \beta \sum_{k_1}\sum_{k_2>k_1}\sum_{k_3>k_2}  \\ 
          & j^2_{k_1 k_2 k_3}(-E_1(\alfa(k_1))E_2(\alfa(k_2))E_1(\alfa(k_3))
              +E_1(\alfa(k_1))E_1(\alfa(k_3)))  \\
          +{}& j^2_{k_2 k_3 k_1}(-E_1(\alfa(k_2))E_2(\alfa(k_3))E_1(\alfa(k_1))
              +E_1(\alfa(k_2))E_1(\alfa(k_1))) \\
          +{}& j^2_{k_3 k_1 k_2}(-E_1(\alfa(k_3))E_2(\alfa(k_1))E_1(\alfa(k_2))
              +E_1(\alfa(k_3))E_1(\alfa(k_2)))\PERIOD
\end{aligned}
\end{equation}
The terms $E_i$ are  defined in  (\ref{e1}-\ref{e4}) and 
the quantities $j^1_{kl}, j^2_{kl}, j^2_{k_1k_2k_3}$ are defined in (\ref{j1}-\ref{j23}).
\item Gibbs measure $\BARIT{\mu}_{M,\beta}^{(2)}(d\eta)
=\frac{1}{\BARIT{Z}_M^{(2)}}
\EXP{-(\BARH_M^{(0)}+\BARH_M^{(1)}+\BARH_M^{(2)})}\BARIT{P}_M(d\eta)$.
\item Relative entropy error
\[
\RELENTR(\BARIT{\mu}_{M,\beta}^{(2)}|\mu_{N,\beta}\circ \COP^{-1})
\sim\BIGO(\epsilon^3)\PERIOD
\]
\end{enumerate}
\end{scheme}

\begin{remark}
{\rm
The estimates in Scheme~\ref{scheme1} and Scheme~\ref{scheme2} refer to 
the equilibrium Gibbs states of lattice systems. Non-equilibrium models such as
the Arrhenius dynamics can be coarse-grained
as shown in \cite{KMVPNAS, KMV}. In this case it is also possible to carry out
a detailed error analysis between   the  exact microscopic   
and the   approximating coarse-grained dynamics. In \cite{kt} an order $\BIGO(\epsilon)$
estimate was proved for the specific relative entropy between microscopic and 
coarse-grained dynamics, while in \cite{KPS} an improved $\BIGO(\epsilon^2)$ estimate was 
shown in the weak topology. Both results are analogous to Scheme~\ref{scheme1}, as the  
invariant Gibbs measure corresponding to the coarse-grained Arrhenius dynamics 
studied there is \VIZ{cg_one}.
}
\end{remark}

We now turn our attention to  the  a posteriori error estimate in Theorem~\ref{mainresult2}. 
The terms derived in Scheme~\ref{scheme2} and the error estimate of 
Theorem~\ref{mainresult} provide an explicit way of calculating the error
made by the 2nd-order coarse-graining.

\smallskip
\noindent{\bf Proof of Theorem~\ref{mainresult2}:}
As we have seen in Proposition~\ref{proposition}, we can write the relative entropy
error as the expectation with respect to the coarse-grained measure of some 
quantities, namely of the difference $\BARH_M-\BARH_M^{(p)}$ for any choice
of $K$.
On the other hand, in the Scheme~\ref{scheme2} we have calculated explicitly
the terms that contribute to the error of $\BARH_M^{(0)}$.
We now conclude the proof  by recalling  the proof of Theorem~\ref{mainresult} and relations 
\VIZ{relative_entropy} and \VIZ{partitions}.
\qed

Using a Metropolis Monte Carlo sampler
(see Appendix~B), we can numerically compute the expectations in \VIZ{relative_entropy}
and \VIZ{partitions} and thus
we can calculate the a posteriori error of 
$\RELENTR(\BARIT{\mu}_{M,\beta}^{(0)}|\mu_{N,\beta}\circ \COP^{-1})$.
Note also that as a consequence the quantity $\EXPECT_{\BARIT{\mu}_{M,\beta}^{(0)}}[R(\eta)]$
can be calculated on-the-fly using the coarse simulation only, at least up to the error
$\BIGO(\epsilon^3)$ according to Scheme~\ref{scheme2}.
\subsection{Coarse-graining  of external fields}
To include the external field in our analysis, we define an effective
external field $\bar{h}$ by
\begin{equation}\label{eff_ext}
\EXP{-\beta\bar{h}(\eta(k))}=\frac{1}{\rho(\eta(k))}\int_{\{F(\sigma|_k)=\eta(k)\}}
\EXP{-\beta\sum_{x\in\CUBE_k}h(x)\sigma(x)}\tilde{\rho}_k(d\sigma)\PERIOD
\end{equation}
for every $k$.
Then the effective full Hamiltonian (for both the interaction and the 
external field contribution) will be given by
\[
\EXP{-\beta\BARH_M(\eta)}=\frac{1}{\rho(\eta)}\int_{\{F(\sigma)=\eta\}}
\EXP{-\beta H_N-\beta\sum_{x\in\Lambda}h(x)\sigma(x)}P_N(d\sigma)
\]
and by perturbing around the coarse-grained Hamiltonian $\BARH_M^{(0)}+\sum_k\bar{h}(\eta(k))$ we obtain
\begin{eqnarray*}
\BARH_M(\eta) & = & \BARH_M{(0)}+\sum_k\bar{h}(\eta(k))\\
&& \mbox{}
-\frac{1}{\beta}
\log\frac{1}{\rho(\eta)}\int_{\{F(\sigma)=\eta\}}\EXP{-\beta\Delta H}
\EXP{-\beta(\sum_k\bar{h}(\eta(k))-\sum_{x\in \Lambda}h(x)\sigma(x))}P_N(d\sigma)\COMMA
\end{eqnarray*}
where $\Delta H=H_N-\BARH_M^{(0)}$.
Since the part that involves the external field can
be written in a product form, we obtain
\begin{eqnarray*}
\BARH_M(\eta)
& = & \BARH_M^{(0)}+\sum_k\bar{h}(\eta(k))-{}\\
&&
{}-\frac{1}{\beta}
\log\int\EXP{-\beta\Delta H}
\prod_k\left(
\EXP{\beta(\bar{h}(\eta(k))-\sum_{x\in \CUBE_k}h(x)\sigma(x))}
\tilde{\rho}_k(d\sigma)\right)\PERIOD
\end{eqnarray*}
The new prior measure is normalized 
due to the definition \VIZ{eff_ext} of the effective external field and therefore 
we can proceed with the cluster expansions as described before.
Eventually, the only issue will be the analytic 
evaluation of the expected values 
$\EXPECT_h[\sigma(x)|\eta]$ with respect to the new prior measure due
to the dependence on the external field.
Furthermore, it is easy to see that in the case of a constant or slowly
varying (at the same scale as the coarse cells) external field
there is no extra contribution to the correction terms of the Hamiltonian.
%
%
%
%
\section{Computational algorithms and numerical experiments}\label{sims}
The 2nd-order approximation for the coarse-grained algorithm as described in
Scheme~\ref{scheme1} has been extensively studied in previous works,
see, e.g., \cite{KMVPNAS,KMV,KPS}, where it has been demonstrated that
it performs well in certain regimes (e.g., long-range interactions).
In this section we present  numerical experiments based on
Scheme~\ref{scheme2} and its comparison
with the 2nd-order coarse-graining Scheme~\ref{scheme1}
for cases where the latter does not give satisfactory results.

Before we present numerical examples we briefly discuss the computational complexity
of the approximations. As a simple measure of complexity we use the number of operations
required for evaluating the Hamiltonian. Although the actual Monte Carlo step does not
require evaluation of the full Hamiltonian the relative complexity with respect to 
the operation count of the full microscopic simulation $q=1$ is properly reflected by
this measure.
Given a potential with the interaction radius $L$ on a $d$-dimensional lattice 
with $N$ sites, the number of operations for evaluation of the microscopic Hamiltonian
$H_N$ is $\BIGO(N L^d)$. The count for the 2nd-order approximation $H^{(0)}_M$ on the coarse lattice
with $M$ sites and coarse-graining ratio $q$ becomes $\BIGO(M L^d/q^d)$. The 3rd-order
approximation $H^{(1)}_M$ involves an additional summation over the interaction range and hence
the operation count is $\BIGO(M L^{2d}/q^{2d})$. In these estimates the ratio $L/q$ is understood
to be equal to one whenever $q\geq L$. In such a case the coarse interactions are reduced to the
nearest-neighbor case. Thus the compression of the interaction kernel $J$ in the corresponding approximations
yields the speed-up of order $\BIGO(q^{2d})$ for $H^{(0)}_M$ and  $\BIGO(q^{3d}/L^d)$ for $H^{(1)}_M$.
We see that the 3rd-order approximation gives an improved error estimate at the same computational
cost whenever $q=L$, in other words, whenever we can compress interactions to the nearest-neighbor 
potential.

We demonstrate the approximation properties in simulations of one-dimensional
Ising-type spin systems. The one-dimensional system provides a suitable
test bed since the exact (analytical) solutions are known for both the classical
Ising system (i.e., nearest-neighbor) and the mean-field model (Curie-Weiss model).
We use the exact solutions to ensure that the simulations are not
influenced by finite-size effects. In all figures the exact solutions visually
coincide with the fully resolved simulations, i.e., $q=1$. We computed
error bars for statistical post-processing, however, they are not displayed in the figures
due to their small relative size as compared to the scales of figures.

In the case of nearest-neighbor interactions the one-dimensional system
does not exhibit phase transition. In fact, the exact solution is given by a well-known
formula (see, e.g., \cite{LavisBell}), which we adopt to our choice of Hamiltonian with
the constant nearest-neighbor ($L=1$) interaction potential of strength $J_0$. 
The equilibrium magnetization curve is then given by
\begin{equation}\label{Isingexact1d}
m_\beta(h) = \frac{\sinh(\beta h)}{\sqrt{\sinh^2(\beta h) 
+ e^{-2 \beta J_0}}}\PERIOD
\end{equation}

On the other hand, for infinitely long attractive interactions there exists a 2nd-order phase
transition and hysteresis behavior is observed according to the global mean-field  
theory for $\beta >\beta_c$, \cite{hm}.
More explicitly, the mean-field (Curie-Weiss) model gives the magnetization curve
as a solution of the non-linear equation
\begin{equation}\label{CWexact1d}
h = \beta J_0 m_\beta - \frac{1}{2\beta}\log\frac{m_\beta}{1-m_\beta}\PERIOD
\end{equation}
The Curie-Weiss model exhibits phase transition at the critical temperature given
by $\beta_c J_0 = 1$ in the case of spins $\{-1,1\}$ ($\beta_c J_0 = 4$ for
spins $\{0,1\}$). 

The approximation of the hysteresis behavior in coarse-grained simulations provides
a good test for derived coarse-graining schemes.
It has been observed previously that hysteresis and critical behavior are not captured properly
for short and intermediate range potentials, \cite{KMV}. 
Similar issues in predicting critical behavior were also observed in \cite{PK06} for
coarse-graining of complex fluids. There an artificial solidification effect was observed for
higher levels of coarse-graining.

In the numerical tests presented
here we demonstrate that the derived corrections improve this behavior even in the
case of nearest-neighbor interactions or high coarse-graining ratio $q$. 
The sampling of the equilibrium measure is done by using microscopic and
coarse-grained Metropolis dynamics discussed briefly in Appendix~B.
We compute isotherms similarly to natural parameter continuation,
i.e., we trace the magnetization $m_\beta$ vs. external field $h$,
first upon increasing the field $h$ from low values 
and then decreasing it from high values.
All simulations have been done with the fine lattice of the size $N=512$.
As derived in Section~\ref{errors} the errors depend on the interplay
of three parameters $q$, $L$ and $\beta$. In examples presented here we 
investigate approximation properties for some regimes in this parameter space.
In the computational examples we choose $J$ to be constant on its support of the
size $2L$, in particular $J(x) =  \frac{J_0}{2L}$, for $|x|\leq L$ and $0$ otherwise.

\begin{figure}
  \centerline{\hbox{\psfig{figure=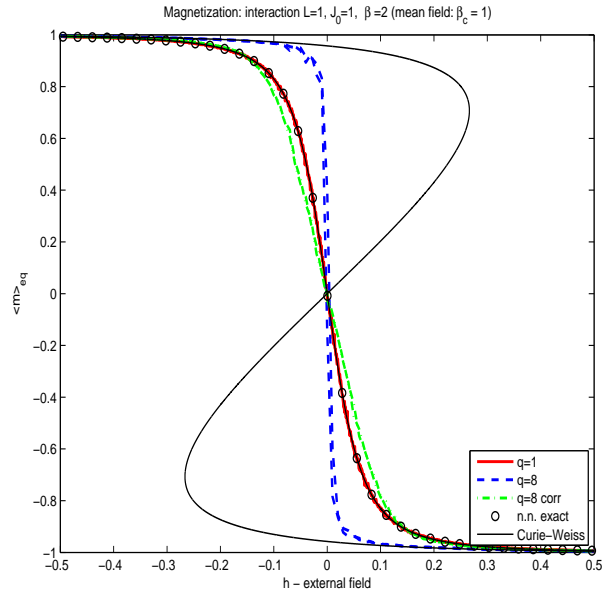,height=8cm,width=8cm}}}
  \caption{\label{L1q8a} Comparison of fully resolved $q=1$ and coarse-grained
           $q=8$ simulations. The interaction range is $L=1$ and the inverse temperature
           is fixed at $\beta=2$.}
\end{figure}

\begin{figure}
  \centerline{\hbox{\psfig{figure=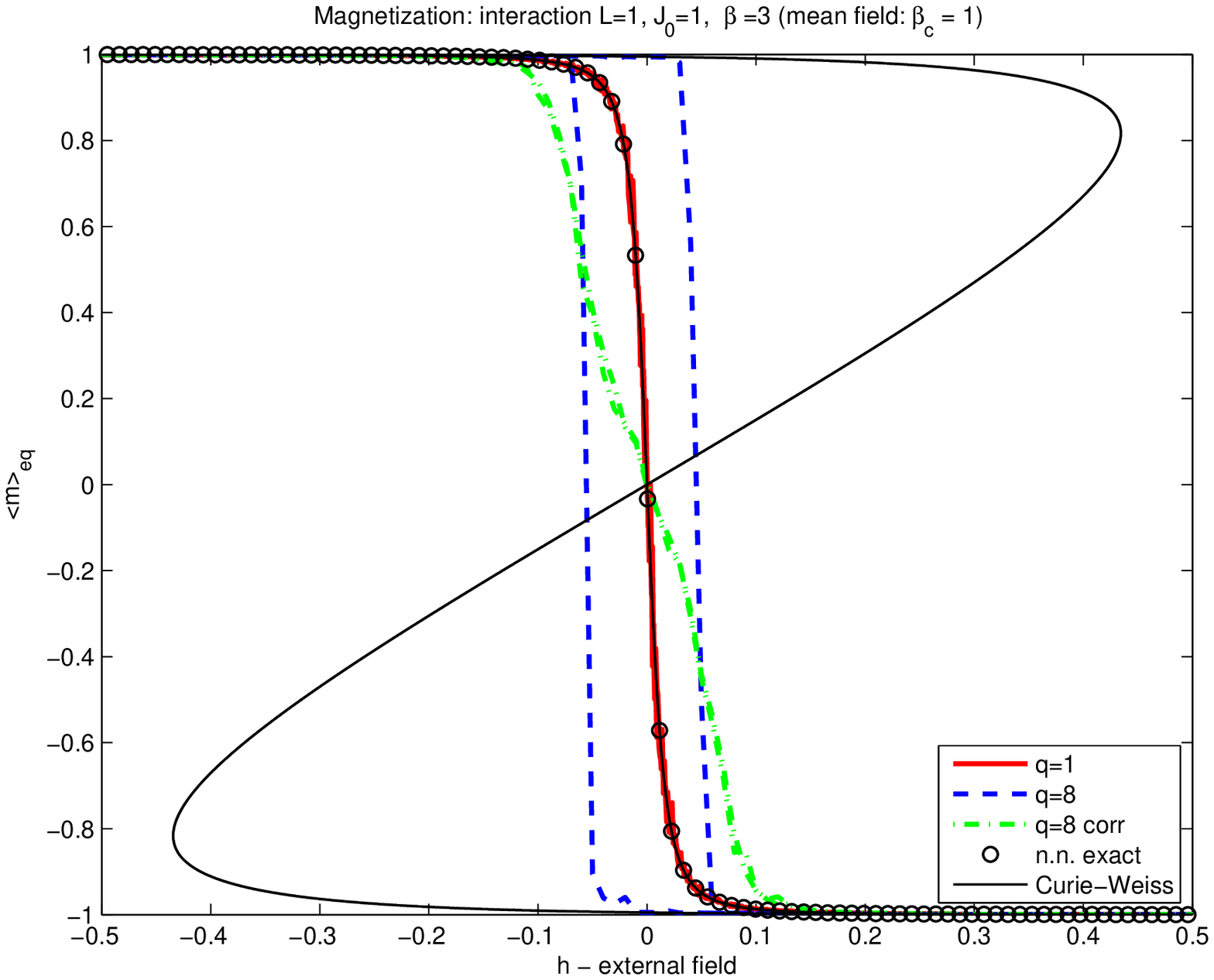,height=8cm,width=8cm}}}
  \caption{\label{L1q8b} Comparison of fully resolved $q=1$ and coarse-grained
           $q=8$ simulations. The interaction range is $L=1$ and the inverse temperature
           is fixed at $\beta=3$.}
\end{figure}

\medskip
\noindent{\it Test Case I: short-range interactions.} We use the classical Ising model with
nearest-neighbor interactions to show efficiency of the 3rd-order correction term.
Figure~\ref{L1q8a}-\ref{L1q8b} depicts simulations at two different temperatures $\beta=2,3$. The
2nd-order approximation exhibits hysteresis behavior for $\beta$ above the critical
value $\beta_c$ for the Curie-Weiss model. It partly follows the
magnetization curve predicted by the mean-field theory. Including the corrections into
the effective Hamiltonian removes the hysteresis and in the case of smaller $\beta$ gives
even a reasonable approximation to the exact solution, which coincides with the case $q=1$
at the scale of the figure.

\begin{figure}
  \centerline{\hbox{\psfig{figure=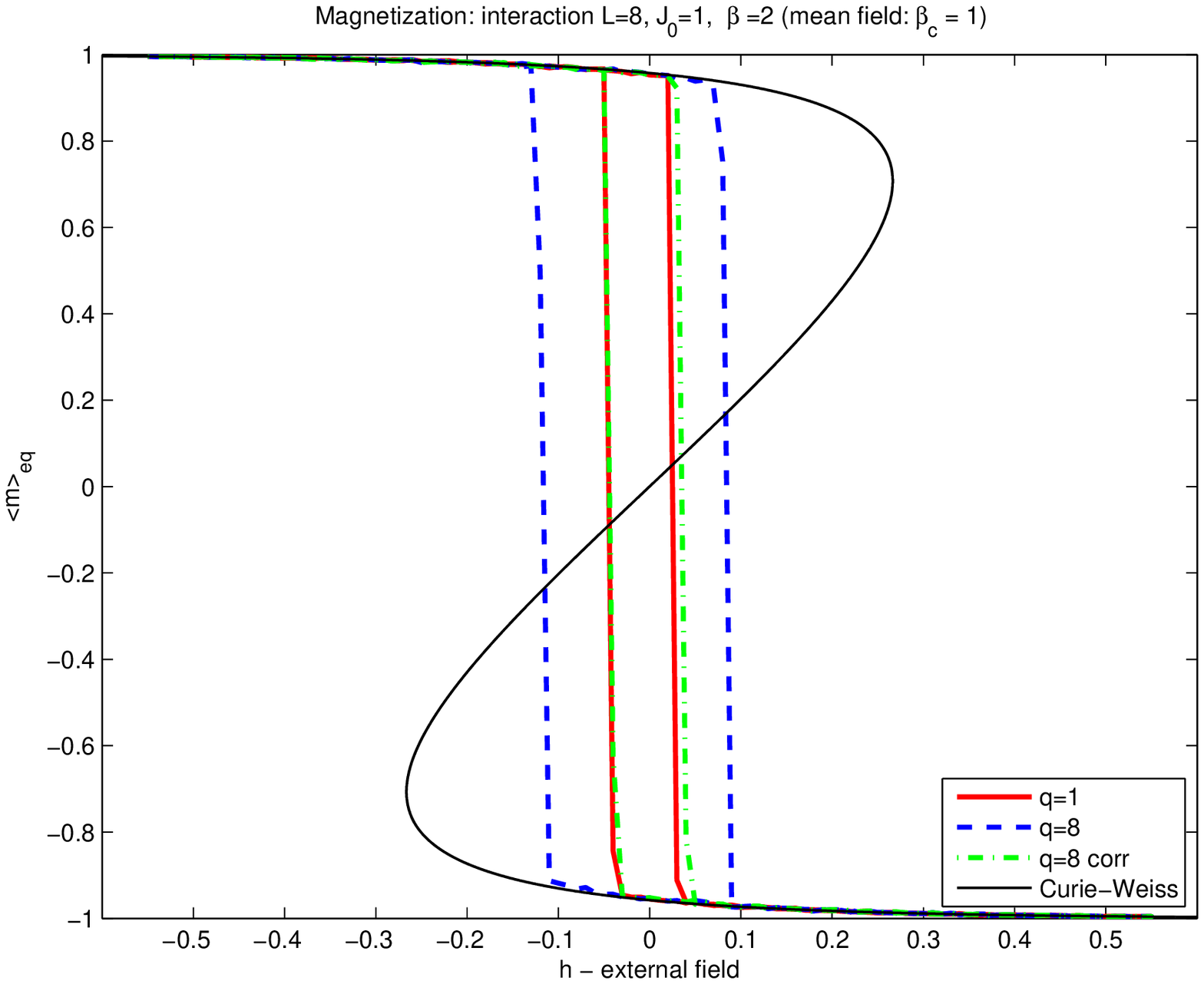,height=8cm,width=8cm}}}
  \caption{\label{L8q8a} Comparison of fully resolved $q=1$ and coarse-grained
           $q=8$  simulations. The interaction range is $L=8$ and the inverse temperature
           is fixed at $\beta=2$.}
\end{figure}

\begin{figure}
  \centerline{\hbox{\psfig{figure=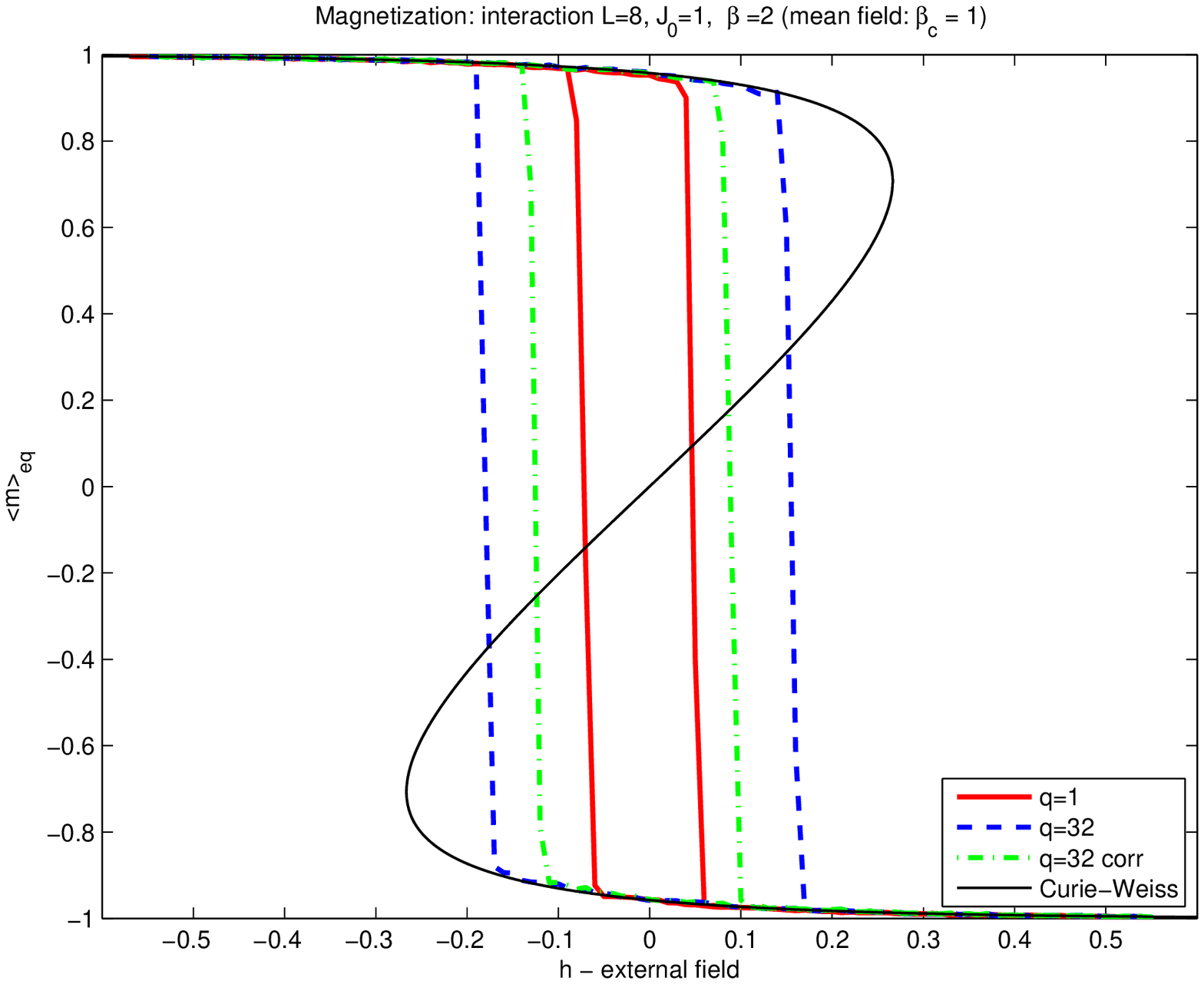,height=8cm,width=8cm}}}
  \caption{\label{L8q8b} Comparison of fully resolved $q=1$ and coarse-grained
           $q=32$ simulations. The interaction range is $L=8$ and the inverse temperature
           is fixed at $\beta=2$.}
\end{figure}

\medskip
\noindent{\it Test Case II: intermediate-range interactions.} In this case
the system exhibits hysteresis behavior (and a 2nd-order phase transition) 
but the correct approximation of the transition between two states is not achieved by
the 2nd-order effective Hamiltonian. In Figure~\ref{L8q8a}-\ref{L8q8b} 
we present an example with extreme coarse-graining
up to the interaction range, i.e., $q=L=8$, and beyond the interaction range $q=32$.
Note that this example is far from the
$\epsilon\sim\frac{q}{L}\ll 1$ limit suggested by Lemma~\ref{lemma1}.
However, as is the case in most asymptotics, computations perform well
even for larger values of $\epsilon$ and especially in this case where
the higher-order corrections are added.

\begin{figure}
  \centerline{\hbox{\psfig{figure=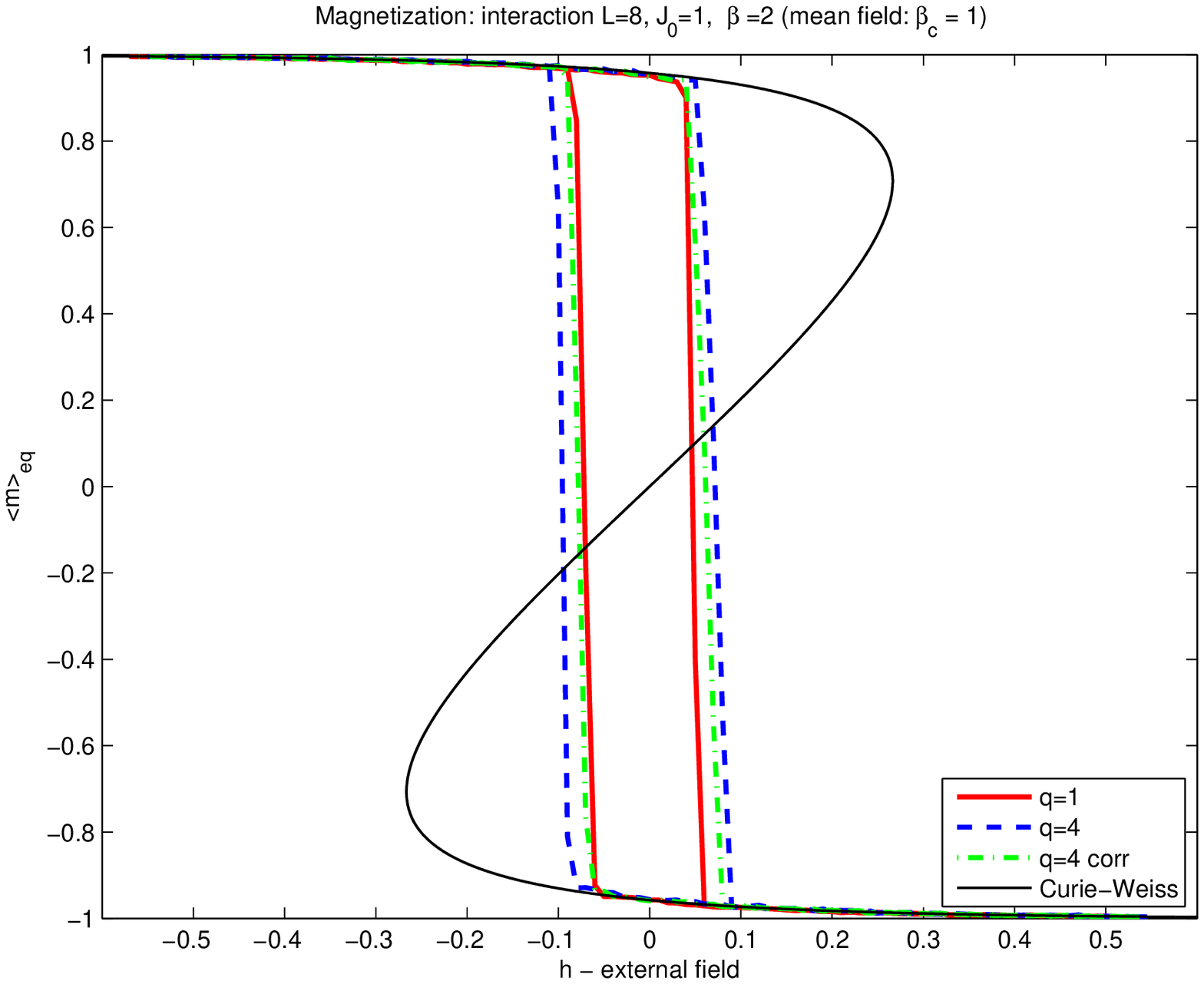,height=8cm,width=8cm}}}
  \caption{\label{L8q4a} Comparison of fully resolved $q=1$ and coarse-grained
           $q=8$ simulations. The interaction range is $L=8$ and the inverse temperature
           is fixed at $\beta=2$.}
\end{figure}

\begin{figure}
  \centerline{\hbox{\psfig{figure=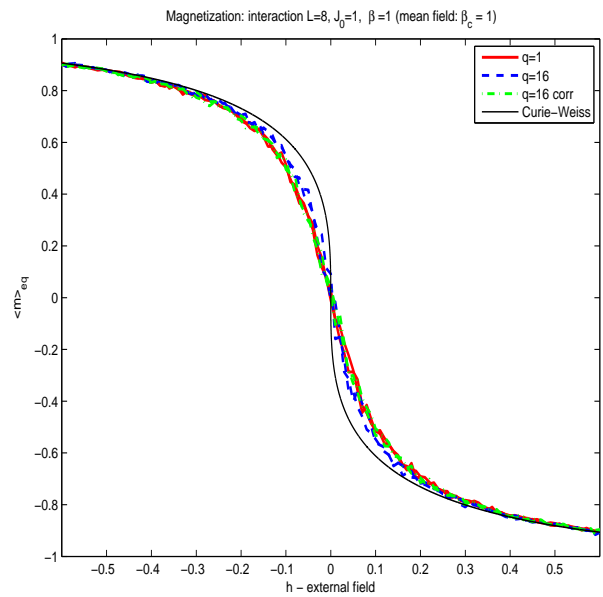,height=8cm,width=8cm}}}
  \caption{\label{L8q4b} Comparison of fully resolved $q=1$ and coarse-grained
           $q=8$ simulations in the ``high temperature'' regime. 
           The interaction range is $L=8$ and the inverse temperature
           is fixed at $\beta=1$.}
\end{figure}

\begin{figure}
  \centerline{\hbox{\psfig{figure=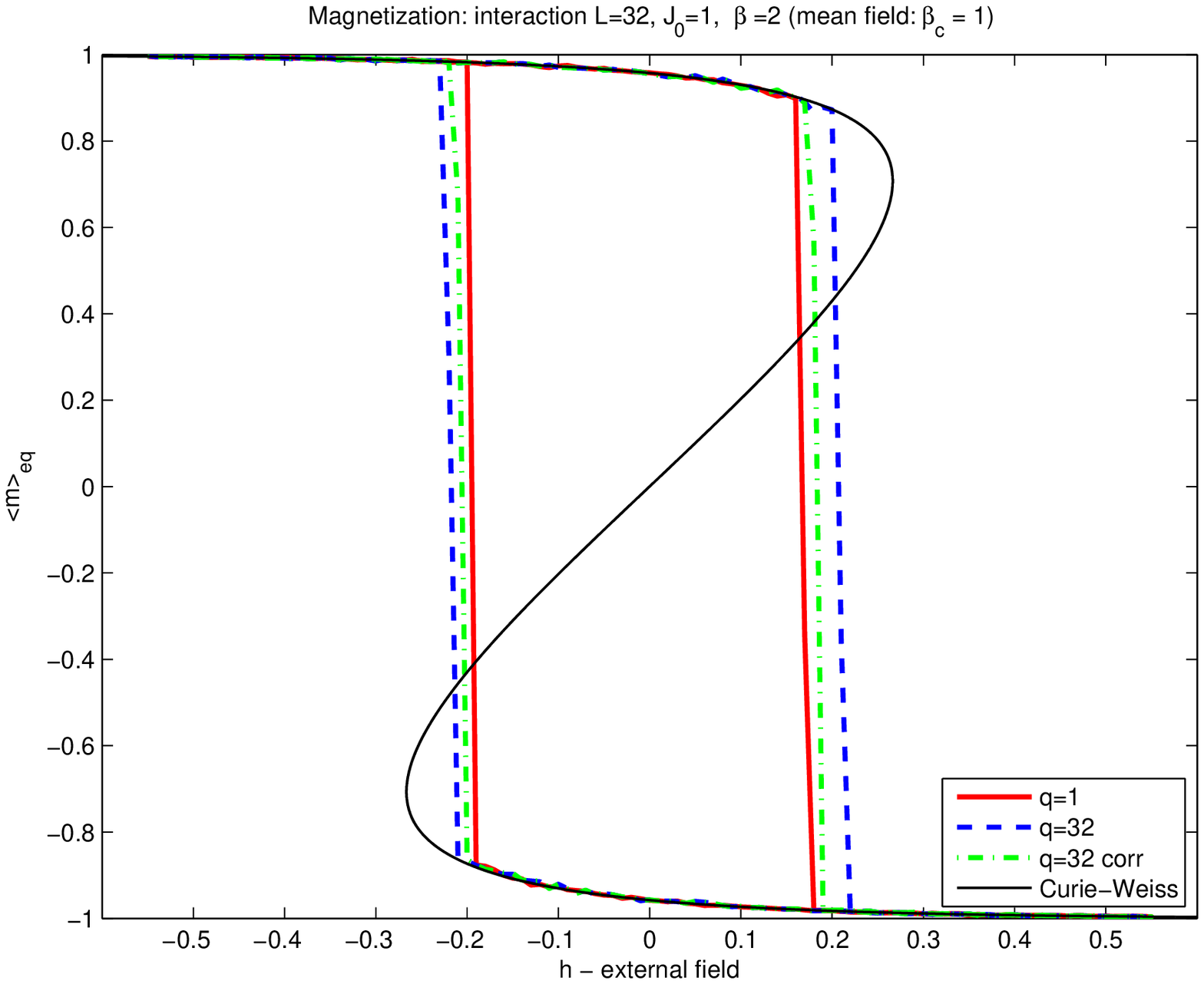,height=8cm,width=8cm}}}
  \caption{\label{L32q32a} Comparison of fully resolved $q=1$ and coarse-grained
           $q=32$  simulations. The interaction range is $L=32$ and the inverse temperature
           is fixed at $\beta=2$.}
\end{figure}

\begin{figure}
  \centerline{\hbox{\psfig{figure=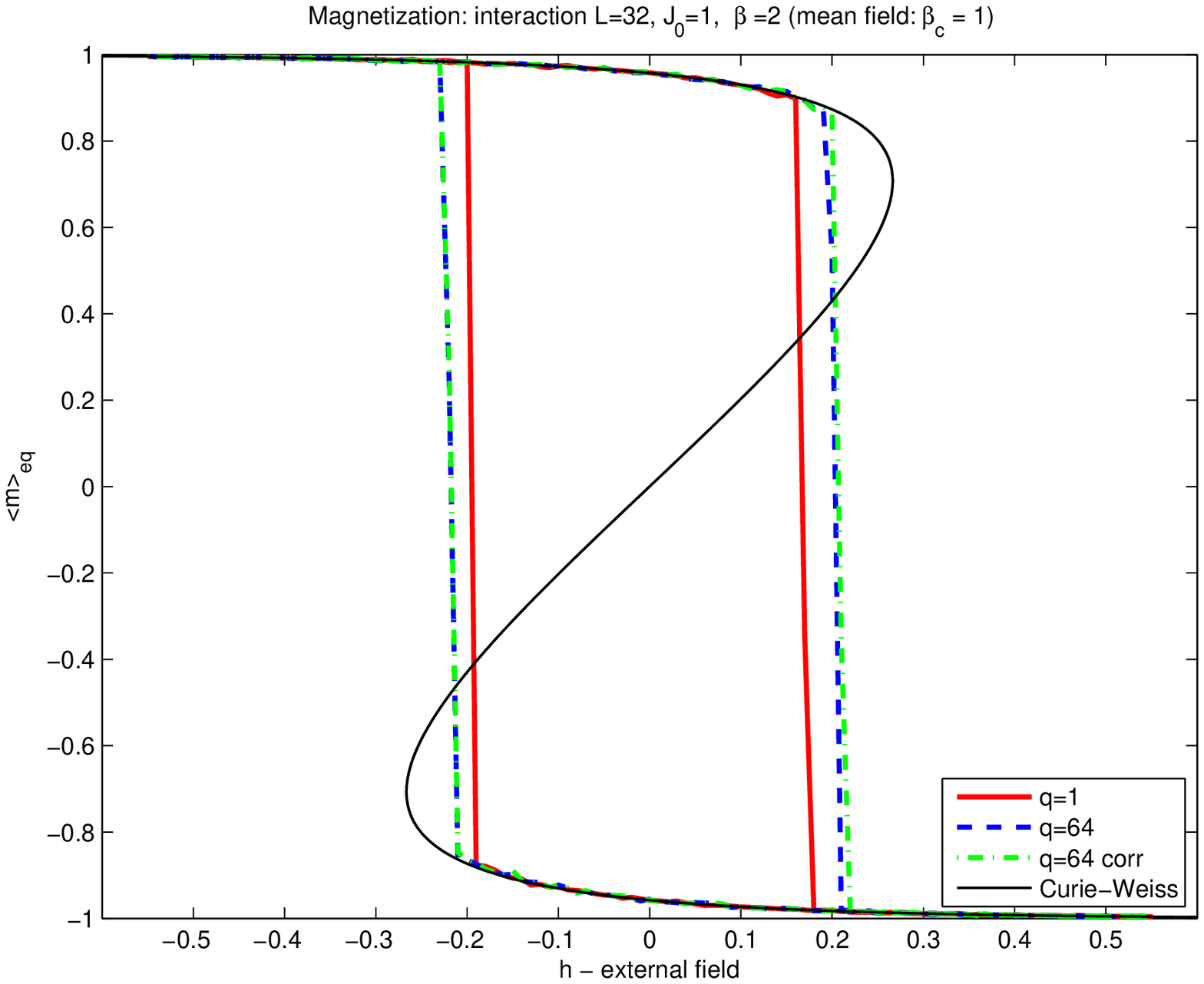,height=8cm,width=8cm}}}
  \caption{\label{L32q32b} Comparison of fully resolved $q=1$ and coarse-grained
           $q=64$ simulations. The interaction range is $L=32$ and the inverse temperature
           is fixed at $\beta=2$.}
\end{figure}

\medskip 
\noindent{\it Test Case III: long-range interactions.} If the interaction range $L$ is sufficiently
large the behavior of the system is well-approximated by the
mean-field solution and coarse-grained simulations will give good predictions already
in the 2nd-order approximation. This is observed in Figure~\ref{L32q32a}-\ref{L32q32b} with coarse-graining
to the range of interactions $q=L=32$. However, coarse-graining beyond the interaction range ($q=64$ in
Figure~\ref{L32q32b})
shows that the 2nd and 3rd-order schemes give almost identical results. 
On the other hand the simulations in Figure~\ref{L8q4a}-\ref{L8q4b} show two regimes where the
2nd-order approximation gives reasonable agreement with the microscopic magnetization
curve.

%
%
%
\newpage
\section*{Appendix A}
In this section we present the detailed calculations involved in obtaining exact formulas 
for the corrections $\bar{H}^{(p)}_M$, $p=1, 2,\dots$ to the coarse-grained
Hamiltonian $\bar{H}^{(0)}_M$.
\SUBSECT{Computation of the term $\int (\Delta_{kk}J(\sigma))^2\tilde{\rho}_k(d\sigma)$}
Using the definition of $\Delta_{kk}J$ we write
\begin{multline*}
\int (\Delta_{kk}J(\sigma))^2\tilde{\rho}_k(d\sigma)= \\
= \sum_{
\genfrac{}{}{0cm}{2}{x,y\in\CUBE_k}{y\neq x}}
\sum_{
\genfrac{}{}{0cm}{2}{x',y'\in\CUBE_k}{y'\neq x'}}
(J(x-y)-\bar{J}(k,k))(J(x'-y')-\bar{J}(k,k))
\EXPECT[\sigma(x)\sigma(y)\sigma(x')\sigma(y')|\alfa]\PERIOD
\end{multline*}
In the above formula the expectation takes different values
according to the following cases
\begin{enumerate}
\item Four particles: $x\neq x',y'$ and $y\neq x',y'$, 
      in which case the expectation gives $E_4(\alfa)$.
\item Three particles: $x=x'$ and $y\neq y'$, or $y=y'$ and $x\neq x'$, or
      $x=y'$ and $x'\neq y$, or $x'=y$ and $x\neq y'$, in which case for some
      distinct $x,y,z$ we obtain 
      $\EXPECT[\sigma^2(x)\sigma(y)\sigma(z)|\alfa]
       =\EXPECT[\sigma(y)\sigma(z)|\alfa]=E_2(\alfa)$, since $\sigma^2(x)=1$.
\item Two particles: $x=x'$ and $y=y'$, or $x=y'$ and $x'=y$, which for the
      distinct particle positions $x,y$ gives
      $\EXPECT[\sigma^2(x)\sigma^2(y)|\alfa]=1$.
\end{enumerate}
We also keep in mind that we already have that $y\neq x$ and $y'\neq x'$.
With the above cases we can substitute for
$\EXPECT[\sigma(x)\sigma(y)\sigma(x')\sigma(y')|\alfa]$
the expression
\begin{eqnarray*}
&& (1-\delta_{x,x'})(1-\delta_{x,y'})(1-\delta_{y,x'})(1-\delta_{y,y'})\times E_4(\alfa)+ \\
&& + \left[\delta_{x,x'}(1-\delta_{y,y'})+\delta_{y,y'}(1-\delta_{x,x'})+\delta_{x,y'}
     (1-\delta_{x',y})+\delta_{x',y}(1-\delta_{x,y'})\right ]\times E_2(\alfa)+ \\
&& + (\delta_{x,x'}\delta_{y,y'}+\delta_{x,y'}\delta_{x',y})\times 1\PERIOD
\end{eqnarray*}
We expand the above products and collect the terms with the factor $\delta_{x,x'}$
(note that the others with the factors $\delta_{x,y'}$ etc. are same by changing 
variables), and
with the factor $\delta_{x,x'}\delta_{y,y'}$ and we obtain
\begin{eqnarray*}
&& 4\sum_{\genfrac{}{}{0cm}{2}{x,y\in\CUBE_k}{y\neq x}}
    \sum_{\genfrac{}{}{0cm}{2}{x',y'\in\CUBE_k}{y'\neq x'}}
        (J(x-y)-\bar{J}(k,k))(J(x'-y')-\bar{J}(k,k))
        \delta_{x,x'}(-E_4(\alfa)+E_2(\alfa))= \\
&& = 4 j^2_{kk}(-E_4(\alfa)+E_2(\alfa))\COMMA\\
&& 2\sum_{\genfrac{}{}{0cm}{2}{x,y\in\CUBE_k}{y\neq x}}
    \sum_{\genfrac{}{}{0cm}{2}{x',y'\in\CUBE_k}{y'\neq x'}}
        (J(x-y)-\bar{J}(k,k))(J(x'-y')-\bar{J}(k,k))
        \delta_{x,x'}\delta_{y,y'}(E_4+1-2E_2)= \\
&& = 2j^1_{kk}(E_4(\alfa)+1-2E_2(\alfa))\COMMA
\end{eqnarray*}

\SUBSECT{Computation of the term $\int (\Delta_{kl}J(\sigma))^2\tilde{\rho}_k(d\sigma)\tilde{\rho}_l(d\sigma)$}
Similarly as in the calculation above we consider all possible cases for the position of 
the particles and we obtain
\begin{eqnarray*}
&& \sum_{\genfrac{}{}{0cm}{2}{x\in\CUBE_k}{y\in\CUBE_l}}
   \sum_{\genfrac{}{}{0cm}{2}{x'\in\CUBE_k}{y'\in\CUBE_l}}
       (J(x-y)-\bar{J}(k,l))(J(x'-y')-\bar{J}(k,l))
       \EXPECT[\sigma(x)\sigma(y)\sigma(x')\sigma(y')|\alfa]=\\
&& 
\begin{aligned}
=\sum_{\genfrac{}{}{0cm}{2}{x\in\CUBE_k}{y\in\CUBE_l}}
    &\sum_{\genfrac{}{}{0cm}{2}{x'\in\CUBE_k}{y'\in\CUBE_l}}
        (J(x-y)-\bar{J}(k,l))(J(x'-y')-\bar{J}(k,l))\times \\
    &\times\Big\{(1-\delta_{x,x'})(1-\delta_{y,y'})\EXPECT[\sigma(x)\sigma(x')|\alfa_k]
         \EXPECT[\sigma(y)\sigma(y')|\alfa_l]+ \\
    &\;\;\;+      \delta_{x,x'}(1-\delta_{y,y'})\EXPECT[\sigma^2(x)|\alfa_k]
          \EXPECT[\sigma(y)\sigma(y')|\alfa_l]+ \\
    &\;\;\;+ \delta_{y,y'}(1-\delta_{x,x'})\EXPECT[\sigma(x)\sigma(x')|\alfa_k]
          \EXPECT[\sigma^2(y)|\alfa_l]+ \\
    &\;\;\;+ \delta_{x,x'}\delta_{y,y'}\EXPECT[\sigma^2(x)|\alfa_k]\EXPECT[\sigma^2(y)|\alfa_l]\Big\}= 
\end{aligned}\\
&& 
\begin{aligned}
= \sum_{\genfrac{}{}{0cm}{2}{x\in\CUBE_k}{y\in\CUBE_l}}
     &\sum_{\genfrac{}{}{0cm}{2}{x'\in\CUBE_k}{y'\in\CUBE_l}}
         (J(x-y)-\bar{J}(k,l))(J(x'-y')-\bar{J}(k,l)) \times \\
     & \times \Big\{\delta_{x,x'}(-E_2(\alfa_k)E_2(\alfa_l)+E_2(\alfa_l))+
          \delta_{y,y'}(-E_2(\alfa_k)E_2(\alfa_l)+E_2(\alfa_k))+\\
     &\;\;\;+  \delta_{x,x'}\delta_{y,y'}(E_2(\alfa_k)E_2(\alfa_l)-E_2(\alfa_l)-E_2(\alfa_k)+1)\Big\}=
\end{aligned}\\
&& = j^2_{kl}(-2E_2(\alfa_k)E_2(\alfa_l)+E_2(\alfa_k)+E_2(\alfa_l)) \\
&&{}+j^1_{kl}(E_2(\alfa_k)E_2(\alfa_l)-E_2(\alfa_l)-E_2(\alfa_k)+1)
\end{eqnarray*}

\SUBSECT{Computation of the term 
$\int \Delta_{kk}J(\sigma)\Delta_{kl}J(\sigma)\tilde{\rho}_k(d\sigma)\tilde{\rho}_l(d\sigma)$}
In the same spirit we have
\begin{eqnarray*}
&& \sum_{\genfrac{}{}{0cm}{2}{x,y\in\CUBE_k}{y\neq x}}
   \sum_{\genfrac{}{}{0cm}{2}{x'\in\CUBE_k}{y'\in\CUBE_l}}
       (J(x-y)-\bar{J}(k,k))(J(x'-y')-\bar{J}(k,l))
       \EXPECT[\sigma(x)\sigma(y)\sigma(x')\sigma(y')|\alfa]=\\
&&
\begin{aligned}
=\sum_{\genfrac{}{}{0cm}{2}{x,y\in\CUBE_k}{y\neq x}}
    &\sum_{\genfrac{}{}{0cm}{2}{x'\in\CUBE_k}{y'\in\CUBE_l}}
     (J(x-y)-\bar{J}(k,k))(J(x'-y')-\bar{J}(k,l)) \times \\
    &\times\Big\{(1-\delta_{x,x'})(1-\delta_{y,x'})\EXPECT[\sigma(x)\sigma(y)\sigma(x')|\alfa_k]
         \EXPECT[\sigma(y')|\alfa_l]+ \\
    &\;\;\;+ \delta_{x,x'}\EXPECT[\sigma^2(x)\sigma(y)|\alfa_k]
         \EXPECT[\sigma(y')|\alfa_l]+
         \delta_{y,x'}\EXPECT[\sigma(x)\sigma^2(y)|\alfa_k]
         \EXPECT[\sigma(y')|\alfa_l]\Big\}=
\end{aligned}\\
&&
\begin{aligned}
=&\sum_{\genfrac{}{}{0cm}{2}{x,y\in\CUBE_k}{y\neq x}}
      \sum_{\genfrac{}{}{0cm}{2}{x'\in\CUBE_k}{y'\in\CUBE_l}}
      (J(x-y)-\bar{J}(k,k))(J(x'-y')-\bar{J}(k,l)) \times\\
    &\times \Big\{-\delta_{x,x'}\LP E_3(\alfa_k)E_1(\alfa_l)-E_1(\alfa_k)E_1(\alfa_l)\RP- \\
    &\;\;\;\;\;\;\;-\delta_{y,x'}\LP E_3(\alfa_k)E_1(\alfa_l)-E_1(\alfa_k)E_1(\alfa_l)\RP\Big\}=
\end{aligned}\\
&& = -2 j^2_{kkl}(E_3(\alfa_k)E_1(\alfa_l)-E_1(\alfa_k)E_1(\alfa_l))\PERIOD
\end{eqnarray*}
Similar computation yields expressions 
$-2j^2_{llk}(E_3(\alfa_l)E_1(\alfa_k)-E_1(\alfa_l)E_1(\alfa_k))$
for the terms $f_{ll}f_{kl}$.

\SUBSECT{Computation of the term 
$\int \Delta_{k_1k_2}J(\sigma)\Delta_{k_2k_3}J(\sigma)\tilde{\rho}_{k_1}(d\sigma)\tilde{\rho}_{k_2}(d\sigma)
\tilde{\rho}_{k_3}(d\sigma)$}
We repeat the same procedure and obtain 
\begin{eqnarray*}
&& \sum_{\genfrac{}{}{0cm}{2}{x\in\CUBE_{k_1}}{y\in \CUBE_{k_2}}}
   \sum_{\genfrac{}{}{0cm}{2}{x'\in\CUBE_{k_2}}{y'\in\CUBE_{k_3}}}
        (J(x-y)-\bar{J}(k_1,k_2))(J(x'-y')-\bar{J}(k_2,k_3))
        \EXPECT[\sigma(x)\sigma(y)\sigma(x')\sigma(y')|\alfa]=\\
&&
\begin{aligned}
  \sum_{\genfrac{}{}{0cm}{2}{x\in\CUBE_{k_1}}{y\in\CUBE_{k_2}}}
     &\sum_{\genfrac{}{}{0cm}{2}{x'\in\CUBE_{k_2}}{y'\in\CUBE_{k_3}}}
           (J(x-y)-\bar{J}(k_1,k_2))(J(x'-y')-\bar{J}(k_2,k_3))\times \\
     &\times \Big\{(1-\delta_{y,x'})\EXPECT[\sigma(x)|\alfa_{k_1}]
             \EXPECT[\sigma(y)\sigma(x')|\alfa_{k_2}]
             \EXPECT[\sigma(y')|\alfa_{k_3}]+\\
     &\;\;\;+\delta_{y,x'}\EXPECT[\sigma(x)|\alfa_{k_1}]\EXPECT[\sigma(y)^2|\alfa_{k_2}]
             \EXPECT[\sigma(y')|\alfa_{k_3}]\Big\}= 
\end{aligned}\\
&& 
\begin{aligned}
  \sum_{\genfrac{}{}{0cm}{2}{x\in\CUBE_{k_1}}{y\in\CUBE_{k_2}}}
     & \sum_{\genfrac{}{}{0cm}{2}{x'\in\CUBE_{k_2}}{y'\in\CUBE_{k_3}}}
           (J(x-y)-\bar{J}(k_1,k_2))(J(x'-y')-\bar{J}(k_2,k_3))\times \\
     &\times \delta_{y,x'}(-E_1(\alfa_{k_1})E_2(\alfa_{k_2})E_1(\alfa_{k_3})+
                         E_1(\alfa_{k_1})E_1(\alfa_{k_3}))=
\end{aligned}\\
&& = j^2_{k_1 k_2 k_3}(-E_1(\alfa_{k_1})E_2(\alfa_{k_2})E_1(\alfa_{k_3})+
     E_1(\alfa_{k_1})E_1(\alfa_{k_3}))\PERIOD
\end{eqnarray*}
%
%
%
\section*{Appendix B: Sampling of the equilibrium measure}
In this section we briefly describe the background of  Monte Carlo
methods  used in Section~\ref{sims} for the sampling of   
the microscopic and coarse-grained equilibrium Gibbs measures.
These algorithms rely on the construction of a suitable Markov Chain  such that its
unique invariant measure is the Gibbs measure we intend to sample. 

\smallskip
\SUBSECT{Microscopic Monte Carlo algorithms}
The sampling algorithm for the microscopic lattice system is given in terms   
a continuous-time jump  Markov process 
that defines a change of the spin  $\sigma(x)$ with the probability $c(x,\sigma)\DT$
over the time interval $[t,t+\DT]$. 
The function $c:\LATT\times \SIGMA\to\R$ is called
a rate of the process. The jump process $\PROCMIC$ is constructed in the following way:
suppose that at the time $t$ the configuration is $\PROCMICRO$, then the probability that
over the time interval $[t,t+\DT]$ the spin at the site $x\in\LATT$ spontaneously 
changes from $\PROCMICRO(x)$ to a new configuration $\sigma^x_{t+\DT}(x)$ is 
$c(x,\sigma)\DT + O(\DT^2)$. We denote the resulting configuration
$\sigma^{x}$.
We require that the dynamics is  such that the invariant measure
of this Markov process is the Gibbs measure \VIZ{microGibbs}. The sufficient condition 
is known as {\it detailed balance} and it imposes a condition on the form
of the rate 
\begin{equation*}
c(x,\sigma) \EXP{-\beta H_N(\sigma)} = 
        c(x,\sigma^{x}) \EXP{-\beta H_N(\sigma^{x})}\PERIOD
\end{equation*}
This condition has a simple interpretation
$c(x,\sigma)$ is the rate of converting $\sigma(x)$ to the value $\sigma^x(x)$
while  $c(x,\sigma^{x})$ is the rate of changing the spin
at the site $x$ back to $\sigma(x)$.
The widely used class of Metropolis-type dynamics satisfies the detailed
balance conditions and has the
rate given by
\begin{equation*}
c(x,\sigma) = G(\beta \Delta_{x}H_N(\sigma))\COMMA\mbox{ where 
        $\Delta_{x}H_N(\sigma) = H_N(\sigma^{x}) - H_N(\sigma)$,}
\end{equation*}
where $G$ is a continuous function satisfying: $G(r)=G(-r)\EXP{-r}$ for
all $r\in\R$. The most common choices in physics simulations are
$G(r)=\frac{1}{1+\EXP r}$ (Glauber dynamics), $G(r)=\EXP{-[r]_+}$, (Metropolis dynamics), with
$[r]_+ = r$ if $r\geq 0$ and $=0$ otherwise, and $G(r)=e^{-r/2}$. Such dynamics are often
used as samplers from the canonical equilibrium Gibbs measure. 

\smallskip
\SUBSECT{Coarse-grained Monte Carlo algorithms}
For the purpose of sampling the coarse-grained Gibbs measure
given by $\BARIT{\mu}_{M,\beta}(\eta)$ we will consider a Markov jump process with 
the rates $ \bar c_b(k,\eta), \bar c_d(k,\eta)$ correspond to the 
addition ("birth")  and removal of a particle ("death") respectively, from the coarse cell 
The resulting process is a {\it birth-death} process $\PROCMAC$ defined 
on the state space $\SPINSPC=\{-q,-q+2, \dots, +q\}$ or $\SPINSPC=\{0,1,\dots q\}$.

As in the case of microscopic dynamics,  we need to ensure the detailed balance condition for the
coarse-grained rates expressed as
\begin{eqnarray*}
&&\bar c_{b}(k,\eta)\BARIT{\mu}_{M,\beta}(\eta)=\bar c_d(k,\eta+\delta_k)\BARIT{\mu}_{M,\beta}(\eta+\delta_k)\\
&&\bar c_d(k,\eta)\BARIT{\mu}_{M,\beta}(\eta)=\bar c_b(k,\eta-\delta_k)\BARIT{\mu}_{M,\beta}(\eta-\delta_k)\PERIOD
\end{eqnarray*}
For the Metropolis-type dynamics
the rates are given by
\begin{eqnarray*}
&&\bar c_{b}(k,\eta) = G(\BARH_M^{(0)}(\eta+\delta_k)-\BARH_M^{(0)}(\eta))
  (q- \eta(k))\COMMA\\
&&\bar c_d(k,\eta) = G(\BARH_M^{(0)}(\eta-\delta_k)-\BARH_M^{(0)}(\eta))
  \eta(k) \PERIOD
\end{eqnarray*}
Notice that with these rates the detailed balance conditions hold
since
\begin{eqnarray*}
&& G(\BARH_M^{(0)}(\eta+\delta_k)-\BARH_M^{(0)}(\eta))
   \EXP{-\BARH_M^{(0)}(\eta)}=G(\BARH_M^{(0)}(\eta-\delta_k)-\BARH_M^{(0)}(\eta))
   \EXP{-\BARH_M^{(0)}(\eta+\delta_k)}\COMMA\\
&& G(\BARH_M^{(0)}(\eta-\delta_k)-\BARH_M^{(0)}(\eta))
   \EXP{-\BARH_M^{(0)}(\eta)}
   =G(\BARH_M^{(0)}(\eta)-\BARH_M^{(0)}(\eta-\delta_k))
   \EXP{-\BARH_M^{(0)}(\eta-\delta_k)}
\end{eqnarray*}
and  $G(r)=G(-r)\EXP{-r}$ for all $r\in\R$.

\smallskip
\SUBSECT{Higher-order coarse-grained Monte Carlo algorithms}
In Theorem~\ref{mainresult} we have suggested higher-order corrections
($\BARH_M^{(p)}$, for $p=1,\ldots$) to the coarse-grained
Hamiltonian $\BARH_M^{(0)}$. Based on this equilibrium theory, and on
the detailed balance condition, we define rates that are capable
of sampling the corrected Gibbs measure corresponding to
$\BARH_M^{(p)}$.
In this case the Metropolis-type rates are given as
\begin{eqnarray*}
&& \bar c_{a}^{(p)}(k,\eta) = G(\BARH_M^{(p)}(\eta+\delta_k)-\BARH_M^{(p)}(\eta))
   (q- \eta(k))\\
&& \bar c_{d}^{(p)}(k,\eta) = G(\BARH_M^{(p)}(\eta-\delta_k)-\BARH_M^{(p)}(\eta))
   \eta(k) \PERIOD
\end{eqnarray*}

%
%
\newpage
%
%

%
%
%
%
%
%

\end{document}